\documentclass{svmult}
\usepackage{amssymb,amsmath,latexsym,amscd}
\usepackage{makeidx}
\usepackage{multicol}

\makeindex

\newcommand{\diag}{\operatorname{diag}}

\def\ga{\mathfrak{a}}

\def\gg{\mathfrak{g}}
\def\gh{\mathfrak{h}}

\def\gk{\mathfrak{k}}

\def\gm{\mathfrak{m}}

\def\gs{\mathfrak{s}}
\def\gt{\mathfrak{t}}

\def\ggg{> \hskip -5 pt >}

\def\C{\mathbb{C}}

\def\H{\mathbb{H}}

\def\R{\mathbb{R}}

\def\cA{\mathcal{A}}

\def\cC{\mathcal{C}}

\def\cI{\mathcal{I}}

\def\cM{\mathcal{M}}

\def\cU{\mathcal{U}}

\newtheorem{thm}[equation]{Theorem}

\newtheorem{lem}[equation]{Lemma}

\newtheorem{defn}[equation]{Definition}

\newtheorem{ex}[equation]{Example}

\newtheorem{rem}[equation]{Remark}

\title*{Infinite Dimensional Multiplicity Free Spaces I: \\ 
Limits of Compact Commutative Spaces}

\institute{Department of Mathematics \\
        University of California \\
        Berkeley, CA 94720--3840, USA \\
        {\tt jawolf@math.berkeley.edu}}

\begin{document}
\author{Joseph A. Wolf} 
\titlerunning{Infinite Dimensional Multiplicity Free Spaces, I}
\maketitle

\begin{abstract}
We study direct limits $(G,K) = \varinjlim\, (G_n,K_n)$ of compact Gelfand 
pairs.  First, we develop a criterion for a direct limit representation 
to be a multiplicity--free discrete direct sum of irreducible representations. 
\index{direct limit}
\index{Gelfand pair}
\index{multiplicity--free}
Then we look at direct limits $G/K = \varinjlim \, G_n/K_n$ of compact 
riemannian symmetric spaces, where we combine our criterion with the 
Cartan--Helgason Theorem to show in general that the regular representation of 
\index{Cartan--Helgason Theorem}
$G = \varinjlim G_n$ on a certain function space $\varinjlim L^2(G_n/K_n)$ 
is multiplicity free.  That method is not applicable for direct limits of
nonsymmetric Gelfand pairs, so we introduce two other methods.  The first,
based on ``parabolic direct limits'' and ``defining representations'', extends
the method used in the symmetric space case.  The second uses some (new) 
branching rules from finite dimensional representation theory.  In both 
cases we define function spaces $\cA(G/K)$, $\cC(G/K)$ and $L^2(G/K)$ to 
which our multiplicity--free criterion applies.  
\end{abstract}

\section{Introduction} \label{sec1}
\setcounter{equation}{0}

Gelfand pairs $(G,K)$, and the corresponding ``commutative'' homogeneous
\index{commutative space}
spaces $G/K$, form a natural extension of the class of riemannian symmetric
spaces.  We recall some of their basic properties.
Let $G$ be a locally compact topological group, $K$ a compact subgroup,
and $M = G/K$.
Then the following conditions are equivalent; see \cite[Theorem 9.8.1]{W2007}.
\medskip

\noindent\phantom{XX}
1. $(G,K)$ is a Gelfand pair, i.e. $L^1(K\backslash G/K)$ is
commutative under convolution.
\hfill\newline\phantom{XX}
2. If $g,g' \in G$ then $\mu_{_{KgK}} * \mu_{_{Kg'K}} =
 \mu_{_{Kg'K}} * \mu_{_{KgK}}$
{\rm (convolution of Dirac measures on $K\backslash G/K$).}
\hfill\newline\phantom{XX}
3. $C_c(K\backslash G/K)$ is commutative under convolution.
\hfill\newline\phantom{XX}
4. The measure algebra $\cM(K\backslash G/K)$ is commutative.
\hfill\newline\phantom{XX}
5. The representation of $G$ on $L^2(M)$ is multiplicity free.
\medskip

\noindent If $G$ is a connected Lie group one can also add
\medskip

\noindent\phantom{XX}
6. The algebra of $G$--invariant differential operators on $M$ is commutative.
\index{invariant differential operator}
\medskip

When we drop the requirement that $K$ be compact, conditions 1, 2, 3 and 4 
lose their meaning because integration on $M$ or $K\backslash G/K$ 
no longer corresponds to integration on $G$.  Condition 5 still makes sense
as long as $K$ is unimodular in $G$.  Condition 6 remains meaningful (and
useful) whenever $G$ is a connected Lie group; there one speaks of
``generalized Gelfand pairs''. 
\medskip

In this note we look at some cases where $G$ and $K$ are not locally
compact, in fact are infinite dimensional, and show in those cases that
the multiplicity--free condition 5 is satisfied.  We first discuss a
multiplicity free criterion that can be viewed as a variation on some
of the combinatoric considerations of \cite{DPW2002}; it emerged from some
discussions with Ivan Penkov in another context.  We then apply the
criterion in the setting of symmetric spaces, proving that direct limits of 
compact symmetric spaces are multiplicity free.  This applies in particular to 
infinite dimensional real, complex and quaternionic Grassmann manifolds,
and it uses some basic symmetric space structure theory.  In particular
our argument for direct limits of compact riemannian symmetric spaces
makes essential of the Cartan--Helgason Theorem, and thus does not
extend to direct limits of nonsymmetric Gelfand pairs.
\medskip

In order to extend the the multiplicity--free result to at least some 
direct limits of nonsymmetric Gelfand pairs, we define the notion of
``defining representation'' for a direct system $\{(G_n,K_n)\}$ where
the $G_n$ are compact Lie groups and the $K_n$ are closed subgroups.
We show how a defining representation for $\{(G_n,K_n)\}$ leads to
a direct system $\{\cA(G_n/K_n)\}$ of
$\C$--valued polynomial function algebras, a
continuous function completion $\{\cC(G_n/K_n)\}$, and a Lebesgue
space completion $\{L^2(G_n/K_n)\}$.  The direct limit spaces
$\cA(G/K)$, $\cC(G/K)$ and $L^2(G/K)$ are the function spaces 
on $G/K = \varinjlim\, G_n/K_n$ which we study as $G$--modules.
\medskip

Next, we prove the multiplicity free property, for the action of
$G$ on $\cA(G/K)$, $\cC(G/K)$ and $L^2(G/K)$, when $\{(G_n,K_n)\}$
is any of six families of Gelfand pairs related to spheres and 
Grassmann manifolds.  Then we go on to prove multiplicity free for
three other types of direct limits Gelfand pairs.  
\medskip

Finally we summarize
the results, extending them slightly by including the possibility 
of enlarging the $K_n$ within their $G_n$--normalizers without losing
the property that $\{K_n\}$ is a direct system.
\medskip

Our proofs of the multiplicity--free condition, for some direct limits of
nonsymmetric Gelfand pairs, use of a number of branching rules, new and old,
for finite dimensional representations.  This lends a certain {\sl ad hoc} 
flavor which I hope can be avoided in the future.
\medskip

Direct limits $(G,K) = \varinjlim (G_n,K_n)$ of riemannian symmetric 
spaces were studied by Ol'shanskii from a very different viewpoint 
\cite{Ol1990}.  He viewed the $G_n$ inside dual reductive pairs and
\index{dual reductive pair}
examined their action on Hilbert spaces of Hermite polynomials.  
Ol'shanskii made extensive use of factor representation theory and 
Gaussian measure, obtaining analytic results on limit--spherical
functions.  See Faraut \cite{Fa2006} for a discussion of spherical functions
in the setting of direct limit pairs.  In contrast to the work of
Ol'shanskii and Faraut we use the rather simple algebraic method of 
renormalizing formal degrees of representations to obtain isometric embeddings 
$L^2(G_n/K_n) \hookrightarrow L^2(G_{n+1}/K_{n+1})$.  
That leads directly to our multiplicity--free results.
\index{formal degree}
\medskip

I am indebted to Ivan Penkov for discussions of multiplicities in direct limit 
representations which are formalized in Theorem \ref{reduction} below.
I also wish to acknowledge hospitality from the Mathematisches
Forschungsinstitut Oberwolfach and support from NSF Grant DMS 04 00420.

\section{Direct Limit Groups and Representations}\label{sec2}

We consider direct limit groups $G = \varinjlim G_n$ and direct limit
representations $\pi = \varinjlim \pi_n$ of them.  This means that
$\pi_n$ is a representation of $G_n$ on a vector space $V_n$, that
the $V_n$ form a direct system, and that $\pi$ is the representation 
of $G$ on $V =\varinjlim V_n$ given by $\pi(g)v = \pi_n(g_n)v_n$ whenever
$n$ is sufficiently large that $V_n \hookrightarrow V$ and $G_n \hookrightarrow
G$ send $v_n$ to $v$ and $g_n$ to $g$.  The formal definition amounts to
saying that $\pi$ is well defined.
\index{direct limit representation}
\medskip

It is clear that a direct limit of irreducible representations is
irreducible, but there are irreducible representations of direct limit
groups that cannot be formulated as direct limits of irreducible
finite dimensional representations.  This is a combinatoric matter and
is discussed extensively in \cite{DPW2002}.  The following definition is
closely related to those combinatorics but applies to a somewhat simpler
situation.

\begin{defn}\label{limit--aligned}{\rm
We say that a representation $\pi$ of $G$ is
{\bf limit--aligned} if it has form $\varinjlim \pi_n$
\index{limit--aligned representation}
in such a way that (i) each $\pi_n$ is a direct sum of primary
representations and (ii) the corresponding
representation spaces $V = \varinjlim V_n$ have the property
every primary subspace of $V_n$ is contained in a primary subspace
of $V_{n+1}$.}
\end{defn}

\begin{thm}\label{reduction}
A limit--aligned representation $\pi = \varinjlim \pi_n$ of
$G = \varinjlim G_n$ is a direct sum of primary representations.
If the $\pi_n$ are multiplicity free then $\pi$ is a multiplicity
free direct sum of irreducible representations.
\end{thm}

\begin{proof}
Let $V = \varinjlim V_n$ be the representation spaces.  Decompose
$V_n = \sum_{\alpha \in I_n} V_{n,\alpha}$ where the $V_{n,\alpha}$
are the subspaces for the primary summands of $\pi_n$.  Write 
$\pi_{n,\alpha}$ for the representation of $G_n$ on $V_{n,\alpha}$,
so $\pi_n = \sum_{\alpha \in I_n} \pi_{n,\alpha}$.
\medskip

Since $\pi$ is limit--aligned, i.e. since each 
$V_{n,\alpha} \subset V_{n+1,\beta}$
for some $\beta \in I_{n+1}$, we may assume $I_n \subset I_{n+1}$ in
such a way that each $V_{n,\alpha} \subset V_{n+1,\alpha}$ for every
$\alpha \in I_n$.  Now $V = \sum_{\alpha \in I} V_\alpha$, discrete 
sum, where $I = \bigcup I_n$ and $V_\alpha = \bigcup V_{n,\alpha}$.  
The sum is direct, for if $u_1 + u_2 + \dots + u_r = 0$ where
$u_i \in V_{\alpha_i}$ for distinct indices $\alpha_1, \dots , \alpha_r$,
then we take $n$ sufficiently large
so that each $u_i \in V_{n,\alpha_i}$ and conclude that 
$u_1 = u_2 = \dots = u_r = 0$.   
Thus $\pi$ is the discrete direct sum of the representations $\pi_\alpha
= \varinjlim \pi_{n,\alpha}$ of $G$ on $V_\alpha$.
\medskip

Let $C_\alpha = \{X: V_\alpha \to V_\alpha \text{ linear } \mid
X\pi_\alpha(g) = \pi_\alpha(g)X \text{ for all } g \in G\}$, the
commuting algebra of $\pi_\alpha$.  If $\pi_\alpha$ fails to be primary 
then $C_\alpha$ contains nontrivial commuting ideals $C_\alpha'$ and
$C_\alpha''$.  Then for $n$ large, the stabilizer $N_{C_\alpha}(V_{n,\alpha})$
of $V_{n,\alpha}$ in $C_\alpha$  contains nontrivial commuting ideals
$N_{C_\alpha'}(V_{n,\alpha})$ and $N_{C_\alpha''}(V_{n,\alpha})$.  That
is impossible because $\pi_{n,\alpha}$ is primary.  We have proved
that $\pi$ is the discrete direct sum of primary representations $\pi_\alpha$.
\medskip

If the $\pi_n$ are multiplicity free then the $\pi_{n,\alpha}$ are
irreducible and it is immediate that the 
$\pi_\alpha = \varinjlim \pi_{n,\alpha}$ are irreducible.
That completes the proof of Theorem \ref{reduction}.
\smartqed\qed\end{proof}

A direct limit of irreducible representations is irreducible, but it is
not immediate that every irreducible direct limit representation can be
rewritten as a direct limit of irreducible representations.  With this
and Theorem \ref{reduction} in mind, we extend Definition \ref{limit--aligned}
as follows.

\begin{defn}\label{lim--irreducible} {\rm
\index{limit--irreducible representation}
A representation $\pi$ of $G = \varinjlim G_n$ is {\bf lim--irreducible}
if it has form $\pi = \varinjlim \pi_n$ where each $\pi_n$ is an irreducible
representation of $G_n$.  Similarly, $\pi$ is {\bf lim--primary} 
if it has form $\pi = \varinjlim \pi_n$ where each $\pi_n$ is a primary
representation of $G_n$.}
\end{defn}

\begin{thm}\label{reduction-cor}
Consider a representation $\pi = \varinjlim \pi_n$ of $G = \varinjlim G_n$
with representation space $V = \varinjlim V_n$.  Suppose that each $\pi_n$
is a multiplicity free direct sum of irreducible highest weight
\index{highest weight representation}
\index{highest weight vector}
representations.  Suppose for $n \ggg 0$ that the direct system map
$V_{n-1} \hookrightarrow V_n$ sends $G_{n-1}$--highest weight vectors 
to $G_n$--highest weight vectors.  Then $\pi$ is a multiplicity free direct
sum of lim--irreducible representations of $G$.
\end{thm}

\begin{proof}  By hypothesis each $\pi_n$ is a direct sum of primary 
representations which, in fact, are irreducible highest weight representations.
We recursively choose highest weight vectors so that
$\pi_{n-1} = \sum \pi_{\lambda,n-1}$ where $\pi_{\lambda,n-1}$ has highest 
weight vector $v_{\lambda,n-1} \in V_{n-1}$ that maps to a highest weight  
vector $v_{\lambda,n} \in V_n$ of an irreducible constituent $\pi_{\lambda,n}$
of $\pi_n$.  This exhibits $\pi$ as a limit--aligned direct sum because it
embeds the summand $V_{\lambda,n-1}$ of $V_{n-1}$ into the irreducible 
summand of $V_n$ that contains $v_{\lambda,n}$. Now
Theorem \ref{reduction} shows that $\pi$ is a multiplicity free direct
sum of lim--irreducible representations of $G$.  
\smartqed\qed\end{proof}

\section{Limit Theorem for Symmetric Spaces}\label{sec3}

We now apply Theorems \ref{reduction} and \ref{reduction-cor} to 
direct limits of compact riemannian symmetric spaces.  
\index{riemannian symmetric space}
Fix a direct system of compact connected Lie groups $G_n$ and subgroups 
$K_n$ such that each $(G_n,K_n)$ is an irreducible riemannian symmetric pair.  
Suppose that the corresponding compact symmetric spaces $M_n = G_n/K_n$ are 
connected and simply connected.  Up to re--numbering and passage to a common
cofinal subsequence the only possibilities are

{\footnotesize
\begin{equation}\label{symmetric-case-class}
\begin{tabular}{|c|l|l|c|c|} \hline
\multicolumn{5}{| c |}
{Compact Irreducible Riemannian Symmetric $M_n = G_n/K_n$} \\
\hline \hline
\multicolumn{1}{|c}{} &
	\multicolumn{1}{c}{$G_n$} &
        \multicolumn{1}{|c}{$K_n$} &
        \multicolumn{1}{|c}{Rank$M_n$} &
        \multicolumn{1}{|c|}{Dim$M_n$} \\ \hline \hline
$1$ & $SU(n)\times SU(n)$ & diagonal $SU(n)$ & $n-1$ & $n^2-1$ \\ \hline
$2$ & $Spin(2n+1)\times Spin(2n+1)$ & diagonal $Spin(2n+1)$ & 
	$n$ & $2n^2+n$ \\ \hline
$3$ & $Spin(2n)\times Spin(2n)$ & diagonal $Spin(2n)$ & 
	$n$ & $2n^2-n$ \\ \hline
$4$ & $Sp(n)\times Sp(n)$ & diagonal $Sp(n)$ & $n$ & $2n^2+n$ \\ \hline 
$5$ & $SU(p+q), \ p = p_n, q = q_n$ & $S(U(p)\times U(q))$ & 
	$\min(p,q)$ & $2pq$ \\ \hline
$6$ & $SU(n)$ & $SO(n)$ & $n-1$ & $\frac{(n-1)(n+2)}{2}$ \\ \hline
$7$ & $SU(2n)$ & $Sp(n)$ & $n-1$ & $2n^2-n-1$  \\ \hline
$8$ & $SO(p+q), \ p = p_n, q = q_n$ & $SO(p) \times SO(q)$ & 
	$\min(p,q)$ & $pq$  \\ \hline
$9$ & $SO(2n)$ & $U(n)$ & $[\frac{n}{2}]$ & $n(n-1)$ \\ \hline
$10$ & $Sp(p+q), \ p = p_n, q = q_n$ & $Sp(p) \times Sp(q)$ & 
	$\min(p,q)$ & $4pq$  \\ \hline
$11$ & $Sp(n)$ & $U(n)$ & $n$ & $n(n+1)$  \\ \hline
\end{tabular}
\end{equation}
}

Fix one of the direct systems $\{(G_n,K_n)\}$ of (\ref{symmetric-case-class}).
Then we have involutive automorphisms
$\theta_n$ of $G_n$ such that the Lie algebras decompose into
$\pm 1$ eigenspaces of the $\theta_n$,
$$
\gg_n = \gk_n + \gs_n \text{ in such a way that } \gk_n = \gg_n \cap \gk_{n+1}
\text{ and } \gs_n = \gg_n \cap \gs_{n+1}.
$$
Then we recursively construct a system of maximal abelian subspaces
$$
\ga_n: \text{ maximal abelian subspace of } \gs_n 
\text{ such that } \ga_n = \gg_n \cap \ga_{n+1}.
$$
The restricted root systems
$$
\Sigma_n = \Sigma_n(\gg_n , \ga_n): 
	\text{ the system of } \ga_n \text{--roots on } \gg_n
$$
form an inverse system of linear functionals: $\Sigma = \Sigma(\gg,\ga)$ 
is the system $\varprojlim \Sigma_n$ of linear functionals on $\ga =
\varinjlim \ga_n$.  In this inverse system, the multiplicities of the
restricted roots will increase without bound, but we can make consistent
choices of positive subsystems
$$
\Sigma^+_n = \Sigma^+_n(\gg_n , \ga_n):
         \text{ system of positive } \ga_n \text{--roots on } \gg_n
$$
so that $\Sigma^+_n \subset \Sigma^+_m|_{\ga_n}$ for $m \geqq n \geqq n_0$.  
Then the corresponding simple root systems
$$
\Psi_n = \Psi_n(\gg_n , \ga_n) = \{\psi_{1,n}, \dots , \psi_{r_n,n}\}:
         \text{ simple } \ga_n \text{--roots on } \gg_n
$$
satisfy $\Psi_n \subset \Psi_m|_{\ga_n}$ for $m \geqq n \geqq n_0$
as well.  Here $r_n = \dim\ga_n$, rank of $M_n$.
\medskip

Recursively define $\theta_n$--stable Cartan subalgebras of 
$\gh_n = \gt_n + \ga_n$ of $\gg_n$  with $\gh_n = \gg_n \cap \gh_{n+1}$.  
Here $\gt_n$ is a Cartan subalgebra of the centralizer $\gm_n$ of $\ga_n$ 
in $\gk_n$.
Now recursively construct positive root systems $\Sigma^+(\gm_n,\gt_n)$
such that if $\alpha \in \Sigma^+(\gm_{n+1},\gt_{n+1})$ then either
$\alpha|_{\gt_n} = 0$ or $\alpha|_{\gt_n} \in \Sigma^+(\gm_n,\gt_n)$.
Then we have positive root systems
$$
\Sigma^+(\gg_n,\gh_n) = \{\alpha \in i\gh_n^* \mid \alpha|_{\ga_n} = 0
\text{ or } \alpha|_{\ga_n} \in \Sigma^+_n(\gg_n , \ga_n) \},
$$
the corresponding simple root systems,
and the resulting systems of fundamental highest weights.
\index{fundamental highest weight}
\medskip

The Cartan--Helgason Theorem says that the irreducible representation
\index{Cartan--Helgason Theorem}
$\pi_\lambda$ of $\gg_n$ of highest weight $\lambda$ gives a summand of
the representation of $G_n$ on $L^2(M_n)$ if and only if (i) 
$\lambda|_{\gt_n} = 0$, so we may view $\lambda$ as an element of
$i\ga_n^*$, and (ii) if $\psi \in \Psi_n(\gg_n , \ga_n)$ then 
$\tfrac{\langle \lambda, \psi \rangle} {\langle \psi, \psi \rangle}$
is an integer $\geqq 0$.  Condition (i) persists under
restriction $\lambda \mapsto \lambda|_{\gh_{n-1}}$ because
$\gt_{n-1} \subset \gt_n$.  Given (i), condition (ii) says that
$\tfrac{1}{2}\lambda$ belongs to the weight lattice of $\gg_n$, so its
restriction to $\gh_{n-1}$ exponentiates to a well defined function on
the corresponding maximal torus of $G_{n-1}$ and thus belongs to the
weight lattice of $\gg_{n-1}$.  Given condition (i) now condition (ii) 
persists under restriction $\lambda \mapsto \lambda|_{\gh_{n-1}}$.  
Define
$$
\Lambda_n = \Lambda(\gg_n,\gk_n,\ga_n) = \left \{\lambda \in i\ga_n^* \left |
	\tfrac{\langle \lambda, \psi \rangle} {\langle \psi, \psi \rangle}
	\text{ integer } \geqq 0 \text{ for all } \psi \in 
	\Psi_n(\gg_n,\ga_n) \right . \right \}.
$$
This is the set of highest weights for representations of $G_n$ on
$L^2(M_n)$, and we have just verified that 
$\Lambda_n|_{\ga_{n-1}} \subset \Lambda_{n-1}$.  Now define the fundamental
highest weights of $\Lambda_n$:
$$
\begin{aligned}
&\xi_{\ell, n} \in i\ga_n^* \text{ defined by } 
	\tfrac{\langle \xi_{\ell, n}, \psi_{m,n} \rangle} 
		{\langle \psi_{m,n}, \psi_{m,n} \rangle} = \delta_{\ell,m}
	\text{ for } 1 \leqq \ell, m \leqq r_n, \\
&\text{ and } \Xi_n = \Xi(\gg_n,\gk_n,\ga_n) = 
	\{\xi_{1,n}, \dots , \xi_{r_n,n}\}.
\end{aligned}
$$
\medskip

\begin{lem}\label{simple-res}
For $n$ sufficiently large, and passing to a cofinal subsequence,
if $\xi \in \Xi_{n-1}$ there is a unique $\xi' \in \Xi_n$ such that
$\xi'|_{a_{n-1}} = \xi$.
\end{lem}

\begin{proof} In the group manifold cases, lines 1, 2, 3 and 4 of
\index{group manifold}
Table \ref{symmetric-case-class}, express $G_n = L_n \times L_n$ and note that
the complexification $(L_{n-1})_\C$ is the semisimple component of a parabolic
subgroup of $(L_n)_\C$.  The restricted root and weight systems of 
$(G_n,K_n)$ are the same as the unrestricted root and weight systems of $L_n$,
and the assertion follows.  
\medskip

In the Grassmann manifold cases, lines 5, 8 and 10 of 
Table \ref{symmetric-case-class}, we first consider the case where 
$\{p_n\}$ is bounded.  Then we may 
assume $p_n = p$ constant and $q_n$ increasing for $n \ggg 0$.
Thus $\ga_{n-1} = \ga_n$, $\Psi_{n-1} = \Psi_n$ (though the multiplicities 
of the restricted roots will increase) and $\Xi_{n-1} = \Xi_n$.  The
assertion now is immediate.
\medskip

In the Grassmann manifold cases we may now assume that both $p_n$ and $q_n$
are unbounded.  If $p_n = q_n$ on a cofinal sequence of indices $n$ we may
assume $p_n = q_n$ for all $n$, so $\Psi_n$ is always of type
$C_{r_n}$.  Then we interpolate pairs and renumber so that 
$p_n = q_n = p_{n-1}+1 = q_{n-1}+1$ for all $n$ and notice that the
Dynkin diagram inclusions $C_{r-1} \subset C_r$ are uniquely determined by
the integer $r$.  If $p_n = q_n$ for only finitely many $n$ and
$p_n < q_n$ on a cofinal sequence of indices $n$ we may assume that
$r_n = p_n < q_n$ for all $n$, so $\Psi_n$ is always of type $B_{r_n}$.
Then we interpolate $(p_{n-1},q{n-1}), (p_{n-1},q_n), (p_{n-1}+1, q_n),
\dots , (p_n,q_n)$ and renumber so that we always have $r_n = r_{n-1}$ or
$r_n = r_{n-1}+1$ and notice that the Dynkin diagram inclusions
$B_{r-1} \subset B_r$ are uniquely determined by the integer $r$.
If $p_n = q_n$ for only finitely many $n$ and also $p_n = q_n$ for only
finitely many $n$ then $p_n > q_n$ on a cofinal sequence of indices $n$,
we may assume $p_n > q_n = r_n$ for all $n$, we interpolate as before
exchanging the r\^ oles of $p_\ell$ and $q_\ell$, and we note again
that the Dynkin diagram inclusions $B_{r-1} \subset B_r$ are uniquely
determined by the integer $r$.  Thus in all cases the fundamental highest 
weights restrict as asserted.
\medskip

In the lower--rank cases, lines 6 and 7 of 
Table \ref{symmetric-case-class}, $\Psi_n$ is of type $A_{n-1}$, so again
restriction to $\ga_{n-1}$ has the required property.  In the hermitian
symmetric case, line 11 of Table \ref{symmetric-case-class}, $\ga_n$ is
a Cartan subalgebra of $\gg_n$ and $\gg_{n-1}$ complexifies to the 
semisimple part of a parabolic subalgebra of $(\gg_n)_\C$, so the
assertion follows as in the group manifold cases.  In the remaining
case, line 9 of Table \ref{symmetric-case-class}, $\Psi_n$ is of type
$C_{n/2}$ for $n$ even, type $B_{(n-1)/2}$ for $n$ odd.  Passing to a cofinal
subsequence we may assume $n$ always even or always odd, and we may 
interpolate as necessary by pairs so that $n$ increases in steps of $2$.
Then, again, there is no choice about the restriction, and the assertion
follows.
\smartqed\qed\end{proof}

In view of Lemma \ref{simple-res}, after passage to a cofinal subsequence
and re-numbering, we may assume the sets $\Xi_n$ ordered so that
\begin{equation}\label{Xi}
\begin{aligned}
\Xi_n = \Xi(\gg_n,\gk_n,\ga_n) =
        &\{\xi_{1,n}, \dots , \xi_{r_n,n}\} \text{ with }\\
&\xi_{\ell, n-1} = \xi_{\ell,n}|_{\ga_{n-1}} \text{ for }
1 \leqq \ell \leqq r_{n-1}.
\end{aligned}
\end{equation}
Now define
\begin{equation}\label{XXi}
\begin{aligned}
&\cI_n: \text{ all } r_n\text{--tuples } I = (i_1,...,i_{r_n})
	\text{ of non--negative integers,}\\
&\cI = \varinjlim \cI_n \text{ where } \cI_n \hookrightarrow \cI_m
	\text{ by } (i_1,...,i_{r_n}) \mapsto (i_1,...,i_{r_n},0,\dots , 0),\\
&\pi_{I,n}: \text{ rep of } G_n \text{ with highest weight }
\xi_I = i_1\xi_1 + \dots + i_p\xi_{r_n}, \\
&\pi_I = \varinjlim \pi_{I,n} \text{ for } I \in \cI.
\end{aligned}
\end{equation}
According to the Cartan--Helgason Theorem, the $\pi_{I,n}$ exhaust 
the representations of $G_n$ on $L^2(M_n)$.  Denote
\begin{equation}\label{XXXi}
V_{I,n}:  \text{ representation space for the abstract representation }
	\pi_{I,n}. 
\end{equation}
Then $V_{I,n}$ occurs with multiplicity $1$ in the
representation of $G_n$ on $L^2(M_n)$.  In effect, 
the representation of $G_n$ on $L^2(M_n)$ is multiplicity free,
and $L^2(M_n) \cong {\bigoplus}_{I\in \cI}  V_{I,n}$ as a $G_n$--module.
However, in the following we must distinguish between 
${\bigoplus}_{I\in \cI}  V_{I,n}$ as a $G_n$--module and $L^2(M_n)$ as
a space of functions.
\medskip

Let $\cU(\gg_n)$ denote the (complex) universal enveloping algebra of
$\gg_n$.  Let $v_{n+1}$ be a highest weight unit vector in $V_{I,n+1}$ for
the action of $G_{n+1}$.  Then we have the $G_n$--submodule
$\cU(\gg_n)(v_{n+1}) \subset V_{I,n+1} \subset L^2(M_{n+1})$.
\medskip

If $u,v \in V_{I,n}$ we write $f_{u,v;I,n}$ for 
$g \mapsto \langle u, \pi_{I,n}(g)v\rangle$, 
the matrix coefficient function on $G_n$.
These matrix coefficient functions span a space $E_{I,n}$ that is
invariant under left and right translations by elements of $G_n$.  As a
$(G_n\times G_n)$--module $E_{I,n} \cong V_{I,n} \boxtimes V_{I,n}^*$.
If $u_n^*$ is the (unique up to scalar multiplication) $K_n$--fixed unit
vector in $V_{I,n}^*$ then the right $K_n$--fixed functions in $E_{I,n}$
form the left $G_n$--module $E_{I,n}^{K_n} \cong V_{I,n}\otimes u_n^*\C
\cong V_{I,n}$.
\medskip

In the following, it is crucial to distinguish between the abstract
representation space $V_{I,n}$ and the space $E_{I,n}^{K_n}$ of functions
on $G_n/K_n$.
\medskip

We normalize Haar measure on $G_n$ (and the resulting measure in $M_n$) to
total mass $1$.  If $u,v, u', v' \in V_{I,n}$ then we have the Schur
Orthogonality relation
$\langle f_{u,v;I,n}, f_{u',v';I,n}\rangle|_{L^2(G_n)}
= (\deg \pi_{I,n})^{-1}\langle u,u'\rangle \overline{\langle v,v'\rangle}$.
\index{Schur Orthogonality}

\begin{thm}\label{fun-restriction}
The space $E_{I,n}^{K_n}$ of functions on $G_n/K_n$ is $G_n$--module
equivalent to $\cU(\gg_n)(v_{n+1} \otimes u_{n+1}^*) 
\subset E_{I,n+1}^{K_{n+1}}$.  We map $E_{I,n}^{K_n}$ into 
$E_{I,n+1}^{K_{n+1}}$ as follows.  \index{formal degree renormalization}
Let $\{w_j\}$ be a basis of $V_{I,n}$ and define
\begin{equation}\label{deg-renormalize}
\begin{aligned}
&\psi'_{n+1,n}\left (\sum c_j\, f_{w_j,u_n^*;I,n}\right ) \\ 
	&\phantom{XXXX} = (\deg \pi_{I,n+1}/\deg \pi_{I,n})^{1/2}
\sum c_j\, f_{w_j,u_{n+1}^*;I,n+1} \in E_{I,n+1}^{K_{n+1}}.
\end{aligned}
\end{equation}  
Then $\psi'_{n+1,n}: E_{I,n}^{K_n} \to E_{I,n+1}^{K_{n+1}}$ is
is $G_n$--equivariant and is isometric for $L^2$ norms on $G_n/K_n$
and $G_{n+1}/K_{n+1}$.  In particular, as $I$ varies with $n$ fixed,
$\psi'_{n+1,n}: L^2(G_n/K_n) \to L^2(G_{n+1}/K_{n+1})$ is a 
$G_n$--equivariant isometry.
\end{thm}

\begin{proof}  We have $a(v_{n+1}) = \xi_I(a)v_{n+1}$ for all 
$a \in \ga$.  The inclusion $G_n \hookrightarrow G_{n+1}$ is 
$G_n$--equivariant, so restriction of functions is $G_n$--equivariant 
and thus is $A$--equivariant, and $(v_{n+1}\otimes u_{n+1}^*)|_{M_n}$ is 
a $\xi_I$--weight vector in $L^2(M_n)$.  If $\alpha$ is a positive 
restricted root for $G_{n+1}$ 
and $e_\alpha \in \gg_{n+1}$ is an $\alpha$ root vector then 
$e_\alpha(v_{n+1}) = 0$.  If $\alpha$ is already a root for $G_n$ and if 
$e_\alpha \in \gg_n$ then we have 
$e_\alpha((v_{n+1}\otimes u_{n+1}^*)|_{M_n}) = 0$.  
Thus either the restriction $(v_{n+1}\otimes u_{n+1}^*)|_{M_n} = 0$ or 
$(v_{n+1}\otimes u_{n+1}^*)|_{M_n}$ is a highest weight vector in 
$E_{I,n}^{K_n}$.
\medskip

Suppose $(v_{n+1}\otimes u_{n+1}^*)|_{M_n} = 0$ as a function on 
$M_n = G_n/K_n$.  Denote $V'_n = \cU(\gg_n)(v_{n+1})$.
It is a cyclic highest weight module for $G_n$ with highest weight $\xi_I$,  
and $(V'_n\otimes u_{n+1}^*\C)|_{M_n} = 0$, and it contains a unique (up to
scalar multiple) $K_n$--invariant unit vector $u'_n$.  
The coefficient function
\index{coefficient function}
$
\varphi(g) := \langle u'_n,\pi_{I,n+1}(g)u'_n \rangle _{V'_n}
        = \int_{G_n} (u'_n\otimes u_n^*)(x) 
		\overline{(u'_n\otimes u_n^*)(x^{-1}g)} dx
$ 
is identically zero because the $u'_n(x)$ factor in the integrand vanishes for 
$x \in G_n$.  But $\varphi|_{G_n}$ is the positive definite 
$(G_n,K_n)$--spherical function on $G_n$ for the representation 
\index{spherical function}
$\pi_{I,n}$, and in particular $\varphi(1) = 1$.  That is a contradiction.  
We conclude that $(v_{n+1}\otimes u_{n+1}^*)|_{M_n} \ne 0$, so 
$(v_{n+1}\otimes u_{n+1}^*)|_{M_n}$ is
a highest weight vector in $E_{I,n}^{K_n}$.  In particular
$E_{I,n}^{K_n} \cong (V'_n\otimes u_{n+1}^*\C)|_{M_n} \subset 
E_{I,n+1}^{K_{n+1}}|_{M_n}$.  That is the
equivariant map assertion.  The unitary map assertion follows by Schur
Orthogonality. \smartqed\qed\end{proof}

Theorem \ref{fun-restriction} gives isometric embeddings
$\psi'_{m,n}: L^2(M_n) \to L^2(M_m)$ for $n \leqq m$.  By construction
$\psi'_{m,n}$ is $G_n$--equivariant.  Define
\begin{equation}
\begin{aligned}
L^2(G/K) = &\varinjlim \{L^2(G_n/K_n), \psi'_{m,n}\}: \text{ direct limit in the}\\
&\text{category of Hilbert spaces and unitary injections.}
\end{aligned}
\end{equation}
We emphasize the renormalizations of Theorem \ref{fun-restriction}.
Without those renormalizations we lose the Hilbert space structure
of $L^2(G/K)$.  

\begin{thm} \label{cor-symm-mfree}
The left regular representation of $G$ on $L^2(G/K)$
is a multiplicity free discrete direct sum of
lim--irreducible representations.  Specifically, that left regular
representation is $\sum_{I \in \cI} \pi_I$ where 
$\pi_I = \varinjlim \pi_{I,n}$ is the irreducible representation
of $G$ with highest weight $\xi_I := \sum i_r\xi_r$.  This
applies to all the the direct systems of {\rm Table \ref{symmetric-case-class}}.
In particular we have the thirteen infinite dimensional multiplicity free spaces
$$
\begin{aligned}
{\rm 1.}\ &SU(\infty) \times SU(\infty)/{\text diag\,}SU(\infty), 
	\text{ group manifold } SU(\infty),\\
{\rm 2.}\ &Spin(\infty)\times Spin(\infty)/{\text diag\,}Spin(\infty), 
	\text{ group manifold } Spin(\infty),\\
{\rm 3.}\ &Sp(\infty)\times Sp(\infty)/{\text diag\,}Sp(\infty), 
	\text{ group manifold } Sp(\infty),\\
{\rm 4.}\ &SU(p + \infty)/S(U(p)\times U(\infty)),\
	\C^p \text{ subspaces of } \C^\infty, \\
{\rm 5.}\ &SU(2\infty)/[S(U(\infty) \times U(\infty))],\
	\C^\infty \text{ subspaces of infinite codim in } \C^\infty, \\
{\rm 6.}\ &SU(\infty)/SO(\infty), \\
{\rm 7.}\ &SU(2\infty)/Sp(\infty), \\
{\rm 8.}\ &SO(p + \infty)/[SO(p)\times SO(\infty)],\text{ oriented } 
	\R^p \text{ subspaces of } \R^\infty, \\
{\rm 9.}\ &SO(2\infty)/[SO(\infty)\times SO(\infty)], \
	\R^\infty \text{ subspaces of infinite codim in } \R^\infty, \\
{\rm 10.}\ &SO(2\infty)/U(\infty), \\
{\rm 11.}\ &Sp(p + \infty)/[Sp(p)\times Sp(\infty)],\
	\H^p \text{ subspaces of } \H^\infty, \\
{\rm 12.}\ &Sp(2\infty)/[Sp(\infty)\times Sp(\infty)], \
	\H^\infty \text{ subspaces of infinite codim in } \H^\infty, \\
{\rm 13.}\ &Sp(\infty)/U(\infty).
\end{aligned}
$$
\end{thm}

\begin{proof}  $\lambda$ is limit--aligned by Theorem 
\ref{fun-restriction}.  Denote $V_I = \bigcup V_{I,n} = \varinjlim V_{I,n}$.  
Then $G$ acts irreducibly on it by $\pi_I = \varinjlim \pi_{I,n}$, and the
various $\pi_I$ are mutually inequivalent because they have different 
highest weights $\xi_I := \sum i_r\xi_r$, and are lim--irreducible
by construction.  Now let 
$V = \sum_{I \in \cI} V_I$.  Then $V = \varinjlim L^2(G_n/K_n) = L^2(G/K)$.
\smartqed\qed\end{proof}

\section{Gelfand Pairs and Defining Representations}

In this section we set the stage for extension of Theorem 
\ref{cor-symm-mfree} to a number of direct systems $\{(G_n,K_n)\}$
of compact nonsymmetric Gelfand pairs.  A glance at \cite{Ya2004} or 
\cite{W2007}
exhibits many such pairs, but here we will only consider those for which
the compact groups $G_n$ are simple.  Here is the Kr\" amer classification 
\index{Kr\" amer classification}
of Gelfand pairs corresponding to compact simple Lie groups 
(see \cite{Kr1979} or \cite{Ya2004} or \cite[Table 12.7.1]{W2007}).

{\small
\begin{equation}\label{kraemer-classn}
\begin{tabular}{|r|l|l|l|l|} \hline
\multicolumn{1}{|l}{} &
\multicolumn{3}{ c |}{$M_n = G_n/H_n$ weakly symmetric} & \multicolumn{1}{ c |}
{$G_n/K_n$ symmetric} \\
\hline \hline
\multicolumn{1}{|l}{} &
\multicolumn{1}{c}{$G_n$} & \multicolumn{1}{|c}{$H_n$}  &
        \multicolumn{1}{|c}{conditions} &
        \multicolumn{1}{|c|}{$K_n$ with $H_n \subset K_n \subset G_n$} \\ \hline \hline
$1$ & $SU(m+n)$ & $SU(m) \times SU(n)]$
        & $n > m \geqq 1$ & $S[U(m) \times U(n)]$ \\
$2$ & $SO(2n)$  & $SU(n)$ & $n$ odd, $n \geqq 3$
        & $U(n)$ \\
$3$ &$E_6$& $Spin(10)$ & & $Spin(10)\cdot Spin(2)$  \\ \hline 
$4$ & $SU(2n+1)$ & $Sp(n)$   & $n \geqq 1$
        & $U(2n) = S[U(2n)\times U(1)]$ \\
$5$ & $SU(2n+1)$ & $Sp(n) \times U(1)$ & $n \geqq 1$
        & $U(2n) = S[U(2n)\times U(1)]$ \\ \hline 
$6$ & $Spin(7)$ & $G_2$ & & (there is none) \\
$7$ & $G_2$ & $SU(3)$ & & (there is none) \\ \hline 
$8$ & $SO(10)$ & $Spin(7) \times SO(2)$ & & $SO(8)\times SO(2)$ \\
$9$ & $SO(9)$ & $Spin(7)$ & & $SO(8)$ \\
$10$ & $Spin(8)$ & $G_2$ & & $Spin(7)$ \\ \hline 
$11$ & $SO(2n+1)$ & $U(n)$ & $n \geqq 2$ & $SO(2n)$ \\
$12$ & $Sp(n)$ & $Sp(n-1) \times U(1)$ & $n \geqq 1$
        & $Sp(n-1) \times Sp(1)$ \\ \hline
\end{tabular}
\end{equation}
}
That gives us the nonsymmetric direct systems $\{(G_n,K_n)\}$ where
\begin{equation} \label{kramer-seq}
\begin{aligned}
&\text{\rm (a) }
G_n = SU(p_n + q_n) \text{ and } K_n = SU(p_n) \times SU(q_n), \ p_n < q_n \\
&\text{\rm (b) }
G_n = SO(2n) \text{ and } K_n = SU(n), \ n \text{ odd, } n \geqq 3 \\
&\text{\rm (c) }
G_n = SU(2n+1) \text{ and } K_n = Sp(n), \ n \geqq 1 \\
&\text{\rm (d) }
G_n = SU(2n+1) \text{ and } K_n = U(1) \times Sp(n), \ n \geqq 1 \\
&\text{\rm (e) }
G_n = SO(2n+1) \text{ and } K_n = U(n), \ n \geqq 2 \\
&\text{\rm (f) }
G_n = Sp(n) \text{ and } K_n = U(1) \times Sp(n-1), \ n \geqq 2
\end{aligned}
\end{equation}

\begin{defn}\label{def-rep}{\rm
Let $\{(G_n,K_n)\}$ be a direct system of Lie groups and closed subgroups.
Suppose that $\pi = \varinjlim \pi_n$ is a lim--irreducible representation
of $G = \varinjlim G_n$, with representation space $V = \varinjlim V_n$,
such that (i) $\pi_n(K_n)$ is the $\pi_n(G_n)$--stabilizer of a vector
$v_n \in V_n$ and (ii) each $v_{n+1} = v_n + w_{n+1}$ where $\pi_n(G_n)$
leaves $w_{n+1}$ fixed.  (Thus the $v_n$ give a coherent system of embeddings
of the $G_n/K_n$.)\,\,  Suppose further that for $n \ggg 0$ the $\pi_n$ have 
the same highest weight vector.  Then we say that $\pi = \varinjlim \pi_n$
is a {\bf defining representation} for $\{(G_n,K_n)\}$.}
\index{defining representation}
\end{defn}

Now let's consider some important examples of defining representations.
We will use those examples later.

\begin{ex}\label{su_suxsu}{\rm
$G_n = SU(p_n + q_n) \text{ and } K_n = SU(p_n) \times SU(q_n)$, $p_n < q_n$,
in (\ref{kramer-seq}).  Let $\{e_1 , \dots , e_{p_n+q_n}\}$ denote the
standard orthonormal basis of $\C^{p_n+q_n}$.  Then $K_n$ is the
$G_n$--stabilizer of $e_1 \wedge \dots \wedge e_{p_n}$ in the representation
$\pi_n = \Lambda^{p_n}(\tau)$ where $\tau$ is the standard (vector)
representation of $SU(p_n + q_n)$ on $\C^{p_n+q_n}$.  In the usual
notation, $e_1 \wedge \dots \wedge e_{p_n}$ also is the highest weight
vector, and the highest weight is $\varepsilon_1 + \dots + \varepsilon_{p_n}$.
If the $p_n$ are bounded, so we may assume each $p_n = p < \infty$, then
$\pi = \varinjlim \pi_n$ is well defined and is a defining representation
for $\{(G_n,K_n)\}$.}\hfill $\diamondsuit$
\end{ex}

\begin{ex}\label{su_sp}{\rm
$G_n = SU(2n+1) \text{ and } K_n = U(1)\times Sp(n), \ n \geqq 1$, in 
(\ref{kramer-seq}).
Again, $\{e_1 , \dots , e_{2n+1}\}$ is the standard orthonormal basis of
$\C^{2n+1}$.  Now $K_n$ is the $G_n$--stabilizer of
$\sum_{\ell = 1}^n e_{2\ell}\wedge e_{2\ell +1}$ in the representation
$\pi_n = \Lambda^2(\tau)$ where $\tau$ is the standard (vector)
representation of $SU(2n+1)$ on $\C^{2n+1}$.  Here $e_1\wedge e_2$ is the
highest weight vector and the highest weight is $\varepsilon_1 +
\varepsilon_2$.  Thus $\pi = \varinjlim \pi_n$ is well defined and 
is a defining representation for $\{(G_n,K_n)\}$.}\hfill $\diamondsuit$
\end{ex}

\begin{ex}\label{so_u}{\rm
$G_n = SO(2n+1) \text{ and } K_n = U(n), \ n \geqq 2$, in (\ref{kramer-seq}).
Let $\{e_1 , \dots , e_{2n+1}\}$ denote the standard orthonormal basis of
$\R^{2n+1}$.  Let $J = \left ( \begin{smallmatrix} 0 & 1 \\ -1 & 0
\end{smallmatrix} \right )$.  Then $K_n$ is the $G_n$--stabilizer of
$\diag\{0, J, \dots , J\} \in \gg_n$ in the adjoint representation of $G_n$,
in other words (in this case) is the $G_n$--stabilizer of
$\sum_{\ell = 1}^n e_{2\ell}\wedge e_{2\ell +1}$ in the representation
$\pi_n = \Lambda^2(\tau)$ where $\tau$ is the standard (vector)
representation of $SO(2n+1)$ on $\R^{2n+1}$.  As in the previous example,
$e_1\wedge e_2$ is the highest weight vector and the highest weight is 
$\varepsilon_1 + \varepsilon_2$.  Thus $\pi = \varinjlim \pi_n$ is well 
defined and is a defining representation for $\{(G_n,K_n)\}$. }
\hfill $\diamondsuit$
\end{ex}

\begin{ex}\label{sp_usp}{\rm
$G_n = Sp(n) \text{ and } K_n = U(1) \times Sp(n-1), \ n \geqq 2$, in
(\ref{kramer-seq}).  In quaternion matrices, $K_n$ is the $G_n$--commutator
of $\diag\{i,0,0,\dots,0\}$.  In $2n \times 2n$ complex matrices it is the 
$G_n$--commutator of $\diag\{J,0, 0, \dots , 0\}$ where 
$J = \left ( \begin{smallmatrix} 0 & 1 \\ -1 & 0 \end{smallmatrix} \right )$.
There, $G_n$ consists of all elements $g \in U(2n)$ such that 
$g\widetilde{J}g^t = \widetilde{J}$, where $\widetilde{J} = 
\diag\{J, J, \dots , J\}$.  Thus $\gg_n$ is given by $x\widetilde{J} +
\widetilde{J} x^t = 0$, and in particular $\diag\{J,0, 0, \dots , 0\} \in \gg_n$.
Now $K_n$ is the $G_n$--stabilizer of $\diag\{J,0, 0, \dots , 0\}$ in the
adjoint representation $\pi_n$ of $G_n$.  That adjoint representation is the
symmetric square of the standard (vector) representation of $G_n$ on $\C^{2n}$,
so it has highest weight $2\varepsilon_1$ and highest weight vector $e_1^2$.
Thus $\pi = \varinjlim \pi_n$ is well defined and is a defining representation 
for $\{(G_n,K_n)\}$. }
\hfill $\diamondsuit$
\end{ex}

\section{Function Algebras} \label{sec5}

Fix a defining representation $\pi = \varinjlim \pi_n$ for 
$\{(G_n,K_n)\}$.  We are going to define algebras
$$
\begin{aligned}
\cA(G_n) &\text{ and } \cA(G) = \bigcup \cA(G_n); \\
&\cA(G_n/K_n) \text{ and } \cA(G/K) = \bigcup \cA(G_n/K_n)
\end{aligned}
$$
\index{function algebra}
of complex--valued polynomial functions and look at their relations to
square--integrable functions.  Let $d_n = \dim_\R V_n$.  Then
we can consider $G_n$ to be a group of real $d_n \times d_n$ matrices.
Since the $G_n$ are reductive linear algebraic groups, this lets us define
\begin{equation}\label{alg}
\begin{aligned}
&\phantom{X} \cA(G_n): \text{ the algebra of all }
        \C\text{--valued functions}\\
         &\phantom{XXXX} f|_{G_n} \text{ where } f: \R^{d_n \times q_n} \to \C
          \text{ is a polynomial,}\\
&\phantom{X} r_n: \cA(G_n) \to \cA(G_{n-1}): \text{ restriction of functions,}\\&\phantom{X} S_n: \text{ kernel of the algebra homomorphism } r_n, \\
&\phantom{X} T_n: G_{n-1}\text{--invariant complement to $S_n$ in } \cA(G_n).
\end{aligned}
\end{equation}
The following is immediate.
\begin{lem} \label{st-decomp}
The restriction $r_n|_{T_n} : T_n \to \cA(G_{n-1})$ is a $G_{n-1}$--equivariant
vector space isomorphism.  In other words we have a $G_{n-1}$--equivariant
injection $(r_n|_{T_n})^{-1}: \cA(G_{n-1}) \hookrightarrow \cA(G_n)$ of
vector spaces with image complementary to the kernel of the restriction
$r_n: \cA(G_n) \to \cA(G_{n-1})$ of functions.  
\end{lem}
Lemma \ref{st-decomp} gives us
$$
\cA(G) = \varinjlim \cA(G_n) = \bigcup \cA(G_n).
$$
Taking the right--invariant functions we arrive at
\begin{equation}\label{d-function-alg}
\begin{aligned}
&\cA(G_n/K_n) := \{h \in \cA(G_n) \mid h(xk) = h(x) \text{ for }
		x \in G_n,\ k \in K_n\}, \\
&\phantom{XX}
\cA(G/K) = \bigcup \cA(G_n/K_n) \\
&\phantom{XX \cA(G/K)} = \{h \in \cA(G) \mid h(xk) = h(x) \text{ for }
		x \in G,\ k \in K\}
\end{aligned}
\end{equation}
These are our basic function algebras.  
\medskip

The algebra $\cA(G_n)$ contains the constants, separates points on $G_n$, and
is stable under complex conjugation.  The Stone Weierstrass Theorem 
\index{Stone Weierstrass Theorem}
is the main component of

\begin{lem}\label{c-density}
The algebra $\cA(G_n)$ is dense in $\cC(G_n)$, the algebra of continuous 
functions $G_n \to \C$ with the topology of uniform convergence.
Let $S'_n$ and $T'_n$ denote the uniform closures of $S_n$ and $T_n$ in
$\cC(G_n)$.  Then $r_n$ extends by continuity to the restriction map
$r'_n: \cC(G_n) \to \cC(G_{n-1})$, that extension $r'_n$ restricts
to a $G_{n-1}$--equivalence $T'_n \cong \cC(G_{n-1})$, 
$\cC(G_n)$ is the vector space direct sum of 
closed $G_{n-1}$--invariant subspaces $S'_n$ and $T'_n$, and
this identifies $\cC(G_{n-1})$ as a $G_{n-1}$--submodule of $\cC(G_n)$.
\end{lem}

\begin{proof}  The density is exactly the Stone--Weierstrass
Theorem in this setting.  Since $S_n$ and $T_n$ involve different sets of
variables, so do $S'_n$ and $T'_n$.  Now $r_n$ extends
to $r'_n$ as asserted and $S'_n$ is the kernel of $r'_n$.  Similarly, 
$S'_n\cap T'_n = 0$, and the induced algebra homomorphism
$r'_n : \cC(G_n) \to \cC(G_{n-1})$ restricts to a $G_{n-1}$--equivariant
map $r'_n : T'_n \cong \cC(G_{n-1})$.  Finally, $S'_n + T'_n$ is 
closed in $\cC(G_n)$ and contains $\cA(G_n)$.  Thus 
$\cC(G_n) = S'_n \oplus T'_n$ and we can identify $\cC(G_{n-1})$ with the
closed $G_{n-1}$--invariant subspace $T'_n$ of $\cC(G_n)$.  
\smartqed\qed\end{proof}

We use the identifications $\cC(G_{n-1}) \subset \cC(G_n)$ of Lemma
\ref{c-density} to form the union $ \bigcup \cC(G_n)$.  Note 
that $ \bigcup \cC(G_n)$ is the algebra of continuous functions on $G$
that depend on only finitely many variables.
Now use the sup norm, thus the topology of uniform convergence, and define
a Banach algebra
$$
\begin{aligned}
\cC(G): &\text{ functions } f:G \to \C \text{ in the uniform limit closure of }
\bigcup \cC(G_n) \\
&\text{ with sup norm and topology of uniform convergence}
\end{aligned}
$$
Passing to the right $K_n$--invariant functions we have Banach function 
algebras
\begin{equation}\label{banach-fn-alg}
\begin{aligned}
&\cC(G_n/K_n) := \{h \in \cC(G_n) \mid h(xk) = h(x) \text{ for }
               x \in G_n,\ k \in K_n \ \text{ and }\\
 &\phantom{XX}
 \cC(G/K) = \bigcup \cC(G_n/K_n) \\
&\phantom{XX\cC(G/K)} = 
	\{h \in \cC(G) \mid h(xk) = h(x) \text{ for } x \in G,\ k \in K\}.
\end{aligned}
\end{equation}
Here $\cA(G_n/K_n)$ is the subalgebra consisting of all $G_n$--finite
functions in $\cC(G_n/K_n)$, and consequently $\cA(G/K)$ is the
subalgebra consisting of all $G$--finite functions in $\cC(G/K)$.  
\medskip

We pass to $L^2$ limits more or less in the same way as (\ref{d-function-alg})
and (\ref{banach-fn-alg}), except that we must rescale to preserve 
$L^2$ norms as in Theorem \ref{fun-restriction}.  For this we need some
machinery from \cite{W2008}.  
Let $\{G_n\}$ be a strict direct system of compact connected Lie groups,
$\{(G_n)_{_\C}\}$ the direct system of their complexifications.
Suppose that, for each $n$,
\begin{equation}
\begin{aligned}
&\text{\,the semisimple part } \ [(\gg_n)_{_\C}, (\gg_n)_{_\C}] \
\text{ of the reductive algebra } \ (\gg_n)_{_\C} \\
&\text{ is the semisimple component of a parabolic
subalgebra of } (\gg_{n+1})_{_\C}.
\end{aligned}
\end{equation}
Then we say that the direct systems $\{G_n\}$ and $\{(G_n)_{_\C}\}$
are {\bf parabolic} and that $\varinjlim G_n$ and $\varinjlim (G_n)_{_\C}$
are {\bf parabolic direct limits}.  This is a special case of the
\index{parabolic direct limit}
\index{parabolic direct system}
definition of parabolic direct limit in \cite{W2005}.
\medskip

Now let $\{G_n\}$ be a strict direct system of compact connected Lie groups
that is parabolic.  We recursively construct Cartan subalgebras
$\gt_n \subset \gg_n$ with $\gt_1 \subset \gt_2 \subset \dots \subset
\gt_n \subset \gt_{n+1} \subset \dots$ and simple root systems
$\Psi_n = \Psi((\gg_n)_{_\C}, (\gt_n)_{_\C})$ such that
each simple root for $(\gg_n)_{_\C}$ is the restriction of exactly one
simple root for $(\gg_{n+1})_{_\C}$.  Then we may assume that
$\Psi_n = \{\psi_{n,1}, \dots ,\psi_{n,p(n)}\}$ in such a way that each
$\psi_{n,j}$ is the $(\gt_n)_{_\C}$--restriction of $\psi_{n+1,j}$
and of no other element of $\Psi_{n+1}$.  The corresponding sets
$\Xi_n = \{\xi_{n,1}, \dots , \xi_{n,p(n)}\}$ of of fundamental highest
weights can be ordered so that they
satisfy: $\xi_{n+1,j}$ is the unique element of $\Xi_{n+1}$ whose
$(\gt_n)_{_\C}$--restriction is $\xi_{n,j}$, for $1 \leqq j \leqq p(n)$.
Exactly as in Theorem \ref{fun-restriction} this gives us isometric
$G_n$--equivariant injections $\psi_{m,n}: L^2(G_n) \to L^2(G_m)$
for $n\leqq m$.  The associated direct limit maps
$\psi_n : L^2(G_n) \to \varinjlim \{L^2(G_n), \psi_{m,n}\}$ define
the direct limit in
the category of Hilbert spaces and unitary maps as the Hilbert space
completion
$$
L^2(G) = \varinjlim_{unitary} \{L^2(G_n), \psi_{m,n}\} =
\left ( \bigcup \psi_n(L^2(G_n)) \right )^{completion}.
$$

\begin{lem}\label{nonsymm-simple-l2}
Let $\{(G_n,K_n)\}$ be one of the systems of {\rm Examples \ref{su_suxsu}, 
\ref{su_sp}, \ref{so_u}} or {\rm \ref{sp_usp}}.  Then $\{G_n\}$ is parabolic
and the $G_n$--equivariant maps
$\psi_{m,n}: L^2(G_n) \hookrightarrow L^2(G_m)$ send 
right--$K_n$--invariants to 
right--$K_{n+1}$--invariants, resulting in $G_n$--equivariant
unitary injections $\psi'_{m,n}: L^2(G_n/K_n) \to L^2(G_m/K_m)$.
\end{lem}

\begin{proof}  We use the defining relations given in Examples 
\ref{su_suxsu}, \ref{su_sp}, \ref{so_u} and \ref{sp_usp}.  In each case 
we look at the subspaces of $L^2$ given by polynomials of degree $\leqq d$; 
those are finite dimensional invariant subspaces of the $\cA(G_n/K_n)$.
We observed above that $\cA(G_n) \hookrightarrow \cA(G_{n+1})$ maps
right--$K_n$--invariants to right--$K_{n+1}$--invariants.  On each 
irreducible summand, the $L^2(G_n) \hookrightarrow L^2(G_{n+1})$ 
differ only by scale from the corresponding summands of $\cA(G_n)$ and
$\cA(G_{n+1})$, so they also map right--$K_n$--invariants to 
right--$K_{n+1}$--invariants. 
\smartqed\qed\end{proof}

\noindent
Now we have some $L^2$ analogs of \ref{d-function-alg} and \ref{banach-fn-alg}.
\begin{equation}\label{l2-fn-alg}
\begin{aligned}
&L^2(G_n/K_n) := \{h \in L^2(G_n) \mid h(xk) = h(x) \text{ for }
               x \in G_n,\ k \in K_n\}, \\
 &\phantom{XX}
 L^2(G/K) = \left ( \bigcup \psi'_n(L^2(G_n/K_n)) \right )^{completion} \\
&\phantom{XXL^2(G/K)} =
        \{h \in L^2(G) \mid h(xk) = h(x) \text{ for } x \in G,\ k \in K\}.
\end{aligned}
\end{equation}
We have $\cA(G/K) \subset \cC(G/K) \subset L^2(G/K)$ for these spaces, and 
$\cA(G/K)$ is the set of polynomial elements in $L^2(G/K)$.

\begin{thm} \label{cor-gf1-mfree}
Let $\{(G_n,K_n)\}$ be one of the direct systems of nonsymmetric Gelfand pairs
given by {\rm Examples \ref{su_suxsu}, \ref{su_sp}, \ref{so_u}} and  
{\rm \ref{sp_usp}}.  Then the left regular representations of 
$G$ on $\cA(G/K)$, $\cC(G/K)$ and $L^2(G/K)$ are multiplicity free discrete 
direct sum of lim--irreducible representations.  In the notation of 
{\rm (\ref{Xi}), (\ref{XXi})} and {\rm (\ref{XXXi})}, those left regular
representations are $\sum_{I \in \cI} \pi_I$ where
$\pi_I = \varinjlim \pi_{I,n}$ is the irreducible representation
of $G$ with highest weight $\xi_I := \sum i_r\xi_r$.  
Thus we have the infinite dimensional multiplicity free spaces
$$
\begin{aligned}
&{\rm (1)} \phantom{X}
SU(p+\infty)/(SU(p)\times SU(\infty)) \text{ for } 1 \leqq p \leqq \infty, \\
&{\rm (2)} \phantom{X} SU(1+2\infty)/(U(1)\times Sp(\infty)), \\
&{\rm (3)} \phantom{X} SO(1 + 2\infty)/U(\infty), \text{ and } \\
&{\rm (4)} \phantom{X} Sp(1+\infty)/(U(1)\times Sp(\infty))
\end{aligned}
$$
\end{thm}

\begin{proof} Examples \ref{su_suxsu}, \ref{su_sp}, \ref{so_u} and 
\ref{sp_usp} have defining representations and well defined function spaces
$\cA(G/K)$ and $\cC(G/K)$.  The same holds for $L^2(G/K)$ by Lemma
\ref{nonsymm-simple-l2}.  In these examples $\{G_n\}$ is parabolic, so
the left regular representations are limit--aligned by 
Theorem \ref{fun-restriction}.  Now the proof of Theorem \ref{cor-symm-mfree}
holds for these four examples, resulting in the multiplicity free property for
their left regular representations.  \smartqed\qed\end{proof}

\section{Pairs Related to Spheres and Grassmann Manifolds} \label{sec6}

In dealing with nonsymmetric Gelfand pairs we have to be very specific
about the embeddings $G_{n-1} \hookrightarrow G_n$, so we 
review a few of those embeddings.
\medskip

\underline{\sc Orthogonal Groups.}
Let $G_n = SO(n_0+2n)$, the special orthogonal group for the 
bilinear form $h(u,v) = \sum_1^{n_0+2n} u_i v_i$.  The embeddings
are given by
$G_n \hookrightarrow G_{n+1}$ given by $x \mapsto \left ( \begin{smallmatrix}
x & 0 & 0\\ 0 & 1 & 0 \\ 0 & 0 & 1 \end{smallmatrix} \right )$.  Then 
$G = \varinjlim G_n$ is the classical direct limit group $SO(\infty)$.
It doesn't matter what $n_0$ is here, but sometimes we have to distinguish
between the cases of even or odd $n_0$, and in any case we want $\{G_n\}$ to
be parabolic, so we jump by two $1$'s instead of just one.
Specifically, this direct system consists either 
of groups of type B (when the $n_0 + 2n$ are odd) or of type D (when the 
$n_0 + 2n$ are even).  In this section $K_n = 
\left \{ \left . \left ( \begin{smallmatrix} 1 & 0 \\ 0 & x
\end{smallmatrix} \right ) \right | x \in SO(n_0+2n-1)\right \} \subset G_n$.
Then $G_n/K_n$ is the sphere $S^{n_0+2n-1}$, 
$G = \varinjlim G_n = SO(\infty)$, and we express
$K = \varinjlim K_n$ as $SO(1)\times SO(\infty - 1)$ to indicate the
embedding $K \hookrightarrow G$.
\medskip

A defining representation for $\{(G_n,K_n)\}$ is given by the family of 
standard (vector) representations $\pi_n$ of $SO(n_0+2n)$ on
$\R^{n_0+2n}$.  Here $\{SO(n_0+2n)\}$ is a parabolic direct system.
In the standard orthonormal basis the $\pi_n$ all have the same  highest 
weight vector $e_1$ and highest weight $\varepsilon_1$.  
Following the considerations of Section \ref{sec5}, this defining 
representation $\pi = \varinjlim \pi_n$ defines the function spaces
$\cA(G_n/K_n)$, $\cC(G_n/K_n)$ and $L^2(G_n/K_n)$.  The $\pi_n$
share a highest weight vector so we have natural equivariant inclusions
$\cA(G_{n-1}/K_{n-1}) \hookrightarrow \cA(G_n/K_n)$,
$\cC(G_{n-1}/K_{n-1}) \hookrightarrow \cC(G_n/K_n)$ and
$L^2(G_{n-1}/K_{n-1}) \hookrightarrow L^2(G_n/K_n)$, and 
thus the limits $\cA(G/K)$, $\cC(G/K)$ and $L^2(G/K)$.  Thus we have the 
regular representation of $G = SO(\infty)$ on those limit spaces.
\medskip

\underline{\sc Unitary Groups.}
Fix $p > 0$ and define $G_n = SU(p+n)$, the special unitary group 
for the complex hermitian form $h(u,v) = 
\sum_1^{p+n} u_i\bar v_i$.  The embedding $G_n \hookrightarrow G_{n+1}$ is 
given by $x \mapsto \left ( \begin{smallmatrix}
x & 0 \\ 0 & 1 \end{smallmatrix} \right )$.  Then $G = \varinjlim G_n$ 
is the classical parabolic direct limit group $SU(\infty)$.  In this section  
$K_n =
        \left \{ \left . \left ( \begin{smallmatrix} 1 & 0 \\ 0 & x
        \end{smallmatrix} \right ) 
	\right | x \in SU(p),\ y \in SU(n) \right \}$.
Then $G_n/K_n$ is a circle bundle over the Grassmann manifold of 
$p$--dimensional linear subspaces of $\C^{p+n}$, 
$G = \varinjlim G_n = SU(\infty)$, and we sometimes express
$K = \varinjlim K_n$ as $SU(p) \times SU(\infty - p)$ to indicate the
embedding $K \hookrightarrow G$.  If $p = 1$ then $G_n/K_n$ is the
sphere $S^{2n+1}$, the complex Grassmann manifold is complex projective
space, and the circle bundle projection is the Hopf fibration.
\vskip .1 cm

Here the defining representation is essentially that of Example \ref{su_suxsu}.
Let $\pi_{\xi_1}$ denote the usual
vector representation of $G_n$ on $\C^{p+n}$.  Write $\pi_{\xi_p}$
for its $p^{th}$ alternating power, representation of $G_n$ on
$\Lambda^p(\C^{p+n})$; it is the first representation of $G_n$ with a
vector fixed under $K_n$.  That
vector is $e_1 \wedge \dots \wedge e_p$ relative to the standard
basis $\{e_1, \dots , e_n\}$ of $\C^n$, and $K_n$ is its $G_n$--stabilizer.
Thus the $\pi_n = \pi_{\xi_p}$ give a defining representation for
$\{(G_n,K_n)\}$.  Note that the $\pi_n$ all have the same highest weight vector 
$e_1 \wedge \dots \wedge e_p$ and highest weight $\varepsilon_1 +
\dots + \varepsilon_p$.
Following the considerations of Section \ref{sec5}, this defining 
representation $\pi = \varinjlim \pi_n$ defines the function spaces
$\cA(G_n/K_n)$, $\cC(G_n/K_n)$ and $L^2(G_n/K_n)$.
The $\pi_n$ share a highest weight vector so we have natural equivariant 
inclusions $\cA(G_{n-1}/K_{n-1}) \hookrightarrow \cA(G_n/K_n)$,
$\cC(G_{n-1}/K_{n-1}) \hookrightarrow \cC(G_n/K_n)$ and
$L^2(G_{n-1}/K_{n-1}) \hookrightarrow L^2(G_n/K_n)$, and
thus the limits $\cA(G/K)$, $\cC(G/K)$ and $L^2(G/K)$.  That gives us 
the regular representation of $G = SU(\infty)$ on those limit spaces.
\medskip

\underline{\sc Symplectic Groups}
Here $Sp(n)$ is the unitary group of the quaternion--hermitian form 
$h(u,v) = \sum_1^n u_i\bar v_i$ on the quaternionic vector space $\H^n$.
We then have $G_n = Sp(n) \times Sp(1)$ where the $Sp(1)$ acts by
quaternion scalars on $\H^n$.  We will also look at its subgroup 
$Sp(n) \times U(1)$ where $U(1)$ is any (they are all conjugate)
circle subgroup of $Sp(1)$, say $\{e^{i\theta} \mid \theta \in \R\}$.
In both cases the embeddings $G_n \hookrightarrow G_{n+1}$ are
specified by the maps $Sp(n) \hookrightarrow Sp(n+1)$ given by
$x \mapsto \left ( \begin{smallmatrix}
x & 0 \\ 0 & 1 \end{smallmatrix} \right )$.  (We are using quaternionic
matrices.)  Then $G = \varinjlim G_n$ is the classical direct limit group 
$Sp(\infty)\times Sp(1)$ and $G' = \varinjlim G'_n$
is $Sp(\infty)\times U(1)$.  (We need the $Sp(1)$ or the $U(1)$ factor
because otherwise, as we will see below, the multiplicity free property 
will fail.)
\medskip

\underline{Symplectic 1.}
First consider the parabolic direct system given by $G_n = Sp(n) \times Sp(1)$.
Given $n$ we have two $Sp(1)$ groups to deal with at the same time,
so we avoid confusion by denoting the $Sp(1)$ factor of $G_n$ as
$Sp(1)_{ext,n}$ ($ext$ for external) and the identity component of the 
centralizer of $Sp(n-1)$ in $Sp(n)$ by $Sp(1)_{int,n}$ 
($int$ for internal).  In our matrix descriptions of $G_n$, the group 
$Sp(1)_{diag,n}$ is the diagonal subgroup in $Sp(1)_{int,n}\times
Sp(1)_{ext,n}$.  Then $G_n = Sp(n) \times Sp(1)_{ext,n}$ and 
$G = \varinjlim G_n = Sp(\infty) \times Sp(1)$.  Now let $K_n = 
	\left \{ \left . \left ( \begin{smallmatrix} 1 & 0 \\ 0 & x
	\end{smallmatrix} \right ) \right | x \in Sp(n-1) \right \} 
	\times Sp(1)_{diag,n}$ and $K = \varinjlim K_n$.
Then $G_n/K_n$ is the sphere $S^{4n-1}$, in other words the Hopf
fibration $3$--sphere bundle over quaternion projective space
$P^{n-1}(\H)$.  In order to indicate the embedding $K \hookrightarrow G$
we express $K$ as $\{1\} \times Sp(\infty - 1) \times Sp(1)$.
\medskip

A defining representation for $\{(G_n,K_n)\}$ is given by the family of
standard (vector) representations $\pi_n$ of $Sp(n)$ on $\C^{2n}$ tensored 
with the standard $2$--dimensional representation of $Sp(1)$ on $\C^2$.  
That representation has an invariant real form $\R^{4n}$. Consider the 
standard orthonormal basis $\{e_i\otimes f_j\}$ of $\C^{2n}\otimes \C^2$.  
The representations $\pi_n$ of $G_n$ there have the same
highest weight vector $e_1\otimes f_1$ and 
highest weight $(\varepsilon_1)_{Sp(n)} + (\varepsilon_1)_{Sp(1)}$.
They give a defining representation for $\{(G_n,K_n)\}$.
Following the considerations of Section \ref{sec5}, this defining 
representation $\pi = \varinjlim \pi_n$ defines the function spaces
$\cA(G_n/K_n)$, $\cC(G_n/K_n)$ and $L^2(G_n/K_n)$.
The $\pi_n$ have the same highest weight vector so we have natural 
equivariant inclusions $\cA(G_{n-1}/K_{n-1}) \hookrightarrow \cA(G_n/K_n)$,
$\cC(G_{n-1}/K_{n-1}) \hookrightarrow \cC(G_n/K_n)$ and
$L^2(G_{n-1}/K_{n-1}) \hookrightarrow L^2(G_n/K_n)$, and
thus the limits $\cA(G/K)$, $\cC(G/K)$ and $L^2(G/K)$.  So we have the regular 
representation of $G = Sp(\infty) \times Sp(1)$ on those limit spaces.
\medskip

\underline{Symplectic 2.}
Next consider the parabolic direct system given by $G_n = Sp(n) \times U(1)$, 
where the $Sp(1)$
factor of $Sp(n) \times Sp(1)$ is replaced by the circle subgroup 
$\{e^{i\theta} \mid \theta \in \R\}$.  Given $n$ we have
two $U(1)$ groups, the $U(1)_{ext,n}$ that is the $U(1)$ factor of
$G_n$ and the corresponding circle subgroup $U(1)_{int,n}$ of
$Sp(1)_{int,n}$.  Then of course we have the diagonal $U(1)_{diag,n}$.
As above we define $K_n$ to be the product group
        $\left \{ \left . \left ( \begin{smallmatrix} 1 & 0 \\ 0 & x
        \end{smallmatrix} \right ) \right | x \in Sp(n-1) \right \} 
	\times U(1)_{diag,n}$ and we set $K = \varinjlim K_n$.
Then $G_n/K_n$ again is the sphere $S^{4n-1}$.  We express
$K$ as $\{1\} \times Sp(\infty - 1) \times U(1)$.
\medskip

A defining representation for $\{(G_n,K_n)\}$ is given by the family of
standard (vector) representations $\pi_n$ of $Sp(n)$ on $\C^{2n}$ tensored
with the standard $1$--dimensional representation of $U(1)$ on $\C$.
The representations $\pi_n$ of $G_n$ there have the same
highest weight vector $e_1\otimes f_1$.  The corresponding
highest weight is $(\varepsilon_1)_{Sp(n)} + (\varepsilon_1)_{U(1)}$, and
the $\pi_n$ give a defining representation for $\{(G_n,K_n)\}$.
Following the considerations of Section \ref{sec5}, this defining 
representation $\pi = \varinjlim \pi_n$ defines the function spaces
$\cA(G_n/K_n)$, $\cC(G_n/K_n)$ and $L^2(G_n/K_n)$.
The $\pi_n$ have the same highest weight vector so we have natural 
equivariant inclusions $\cA(G_{n-1}/K_{n-1}) \hookrightarrow \cA(G_n/K_n)$,
$\cC(G_{n-1}/K_{n-1}) \hookrightarrow \cC(G_n/K_n)$ and
$L^2(G_{n-1}/K_{n-1}) \hookrightarrow L^2(G_n/K_n)$, and
thus the limits $\cA(G/K)$, $\cC(G/K)$ and $L^2(G/K)$.  So we have the regular 
representation of $G = Sp(\infty) \times U(1)$ on those limit spaces.
\medskip

\underline{Symplectic 3.}
A variation on the case just considered is where $K_n =
	\left \{ \left . \left ( \begin{smallmatrix} z & 0 \\ 0 & x
	\end{smallmatrix} \right ) \right | z \in U(1), x \in Sp(n-1) \right \}
	\times U(1)$ and $K = \varinjlim K_n$.
Then the $U(1)$ factor of $G_n$ is contained in $K_n$ so it acts trivially
on $G_n/K_n$.  Thus $G_n/K_n$ is a $2$--sphere bundle over $P^{n-1}(\H)$
exactly as in the ``Symplectic 2'' case.  We express
$K$ as $U(1)\times Sp(\infty - 1) \times U(1)$.
The groups $K_n$ are larger than the case ``Symplectic 2'' just considered,
so the present function spaces $\cA(G_n/K_n)$, $\cC(G_n/K_n)$ and 
$L^2(G_n/K_n)$ are subspaces of those of ``Symplectic 2'', and the
same holds for their limits $\cA(G/K)$, $\cC(G/K)$ and $L^2(G/K)$.
Now we have the regular representation of $G = Sp(\infty) \times Sp(1)$
on those limit spaces.
\medskip

\underline{Symplectic 4.} 
A variation on the ``Symplectic 1'' case is where 
$G_n = Sp(n) \times Sp(1)$ and $K_n =
        \left \{ \left . \left ( \begin{smallmatrix} z & 0 \\ 0 & x
        \end{smallmatrix} \right ) \right | z \in U(1), x \in Sp(n-1) \right \}
        \times Sp(1)$ and $K = \varinjlim K_n$.
Then the $Sp(1)$ factor of $G_n$ is contained in $K_n$, so it acts
trivially on $G_n/K_n$.  Thus $G_n/K_n = Sp(n)/[U(1)\times Sp(n-1)]$ 
is a $2$--sphere bundle over $P^{n-1}(\H)$, exactly as in the
``Symplectic 3'' case above.  We express
$K$ as $U(1)\times Sp(\infty - 1) \times Sp(1)$ and we note that the
function spaces $\cA(G_n/K_n)$, $\cC(G_n/K_n)$ and $L^2(G_n/K_n)$
are exactly the same as those of ``Symplectic 3'', so the same 
holds for their limits $\cA(G/K)$, $\cC(G/K)$ and $L^2(G/K)$.
Thus we have the regular representation of $G = Sp(\infty) \times Sp(1)$
on those limit spaces.
\medskip

The classifications of Kr\" amer
\cite{Kr1979} and Yakimova \cite{Ya2004} (see \cite{W2007}) show that the six 
direct systems just described, one orthogonal, one unitary, and four 
symplectic, all consist of Gelfand pairs.
\medskip

\section{Limits Related to  Spheres and Grassmann Manifolds}\label{sec7}

In this section we prove a the multiplicity free property for the
direct limits of Gelfand pairs described in Section \ref{sec6}.
\index{Grassmann manifold}

\begin{thm}\label{spec}
Let $(G,K) = \varinjlim \{(G_n, K_n)\}$ where $\{(G_n,K_n)\}$ is one
of the six systems described in {\rm Section \ref{sec6}}.  Let $\cA(G/K)$, 
$\cC(G/K)$ and $L^2(G/K)$ be as described there.  Then the regular
representations of $G$ on $\cA(G/K)$, $\cC(G/K)$ and $L^2(G/K)$ are
multiplicity--free discrete direct sums of lim--irreducible representations.
\end{thm}

\begin{proof}
We run through the proof of Theorem \ref{spec} for the three types of
limit groups $G$.  In each case we do this by examining the representation
of $G_n$ on $\cA(G_n/K_n)$, verifying the limit--aligned condition, and
applying Theorem \ref{reduction-cor} to the regular representation of
$G$ on $\cA(G/K)$.  We already know the result for the orthogonal group
case, where the $(G_n,K_n)$ are symmetric pairs, but we need the
representation--theoretic information from that case in order to deal
with the other cases.
\medskip
 
{\bf Orthogonal Group Case.}
Here we shift the index so that $G_n = SO(n)$ and $K_n = SO(n-1)$. 
Then $G_n/K_n$ is the unit sphere in $\R^n$.  The $G_n$--finite functions 
on $G_n/K_n$ are just the restrictions of polynomial functions on $\R^n$.  
Let $\psi_{1;n}$ denote the usual representation of $G_n$ on $\R^n$ and let
$\xi$ denote its highest weight.  Choose orthonormal linear coordinates 
$\{x_1, \dots , x_n\}$ of that $\R^n$ such that the monomial $x_1$ is a
highest weight vector.  Then the representation of $G_n$ on the
space of polynomials of pure degree $\ell$ is of the form 
$\psi_{\ell;n} \oplus \gamma_{\ell;n}$ where $\psi_{\ell;n}$ is the 
irreducible representation of highest weight $\ell\xi$ and highest 
weight vector $x_1^\ell$. Then $\gamma_{\ell;n}$ is 
the sum of the $\psi_{\ell - 2j;n}$ for $1 \leqq j \leqq [\ell/2]$, 
and the representation space of that $\psi_{\ell - 2j;n}$ consists of 
the polynomial functions on $\R^n$ divisible by $||x||^{2j}$ but not 
by $||x||^{2j+2}$.  Write $E_{\ell;n}$ for the space of functions on 
$G_n/K_n$ obtained by restricting those polynomials of degree $\ell$ 
contained in the representation space for $\psi_{\ell;n}$.  
Then $\cA(G_n/K_n) = \sum_{\ell \geqq 0}\, E_{\ell;n}$.
\medskip

We now verify that the inclusions 
$\cA(G_n/K_n) \hookrightarrow \cA(G_{n+1}/K_{n+1})$
send $E_{\ell;n}$ into $E_{\ell;n+1}$, so that the representation of $G$
on $\cA(G/K)$ is limit--aligned and Theorem \ref{reduction-cor} shows
that $\varinjlim \cA(G_n/K_n)$ is the multiplicity free direct sum of 
lim--irreducible $G$--modules $E_\ell = \varinjlim E_{\ell;n}$.  
For that, note that
the restriction $\cA(G_{n+1}/K_{n+1}) \to \cA(G_n/K_n)$ is obtained by
setting $x_{n+1}$ equal to zero.  Thus the inclusions
$E_{\ell,n} \hookrightarrow \cA(G_{n+1}/K_{n+1})$ are given by 
identifying the function $x_1^\ell : \R^n \to \R$ with the function
$x_1^\ell : \R^{n+1} \to \R$ and applying $G_n$--equivariance.  Now
$\cA(G/K) = \varinjlim \cA(G_n/K_n)$
is the direct sum of the $E_\ell = \varinjlim E_{\ell;n}$, and the
representations of $G$ on the $E_\ell$ are the mutually inequivalent
lim--irreducible $\varinjlim \psi_{\ell;n}$.  That gives an elementary
proof for the case $G = SO(\infty)$ and $K = SO(\infty - 1)$.
\medskip

{\bf Unitary Group Cases.} 
Here we shift the index so that $G_n = SU(n)$ and $K_n = SU(p)\times SU(n-p)$,
$n > p$.  So $G = SU(\infty)$ and $K = SU(p) \times SU(\infty - p)$.
Without loss of generality assume $n > 2p$ so that the $(G_n,K_n)$ are
Gelfand pairs.  Recall the defining representation $\pi = \varinjlim \pi_n$
where $\pi_n = \pi_{\xi_p}$, the $p^{th}$ exterior power of the vector 
representation of $G_n$ on $\C^n$.  So $K_n$ is the $G_n$--stabilizer of
$e_{I_0}:= e_1\wedge \dots \wedge e_p$, resulting in the map
$G_n/K_n \hookrightarrow \Lambda^p(\C^n)$ by $gK_n \mapsto g(e_{I_0})$.
\medskip

We have $\C$--linear functions $z_I$ 
on $\Lambda^p(\C^n)$ dual to the basis of $\Lambda^p(\C^n)$ consisting of
the $e_I$ with $I = (i_1, \dots , i_p)$ where $1 \leqq i_1 < \dots < i_p \leqq
n$.  (Here $I_0 = (1,2,\dots,p)$.)
Their real and imaginary parts generate the algebra $\cA(G_n/K_n)$.
Relative to the diagonal Cartan subalgebra of $\gg_n$ the $e_I$ are
weight vectors, and $e_{I_0}$ is the highest weight vector, for $\pi_{\xi_p}$.
Now the action of $G_n$ on the polynomials of degree $\ell$ in the
$z_I$ and the $\overline{z_I}$ is 
$\sum_{r+s = \ell}\pi_{r\xi_p + s\xi_{n-p}}$, where $\pi_{r\xi_p + s\xi_{n-p}}$
has highest weight $r\xi_p + s\xi_{n-p}$ and highest weight vector
$\overline{z_{I_0}^r} z_{I_0}^s$.  Those representations are mutually
inequivalent, using $n > 2p$, and $\cA(G_n/K_n) = \sum_{\ell \geqq 0}
\sum_{r+s = \ell} E_{r,s;n}$ where $G_n$ acts on $E_{r,s;n}$ by
$\pi_{r\xi_p + s\xi_{n-p}}$.  The $\cA(G_n/K_n) \hookrightarrow
\cA(G_{n+1}/K_{n+1})$ are given on the level of $E_{r,s;n} \hookrightarrow
E_{r,s;n+1}$ by identifying 
$\overline{z_{I_0}^r} z_{I_0}^s: \Lambda^p(\C^n) \to \C$
with 
$\overline{z_{I_0}^r} z_{I_0}^s: \Lambda^p(\C^{n+1}) \to \C$.
In view of Theorem \ref{reduction-cor},
it follows that the representation of $G$ on $\cA(G/K)$ is a 
limit aligned discrete direct sum of mutually inequivalent lim--irreducible 
representations.
\medskip

We will need the case $p = 1$ when we look at the symplectic group cases.
There $G_n = SU(n)$ and $K_n = \{1\} \times SU(n-1)$, and the $G_n$--finite
functions on $G_n/K_n$ are just the restrictions of finite linear
combinations of the functions $z^r\bar z^s$.  We saw how to
decompose $\cA(S^{2n-1})$ into irreducible modules for $SO(2n)$: it
is the sum of the spaces $E_{\ell;2n}$ described above with highest weight
$\ell\xi$ and highest weight vector $x_1^\ell$, where, of course,
$x_j = \frac{1}{2}(z_j + \bar z_j)$.  In terms of the Dynkin diagram 
that representation is 

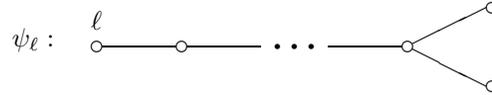
\centerline{
\setlength{\unitlength}{.75 mm}
\begin{picture}(90,20)
\put(0,10){$\psi_\ell: $}
\put(15,10){\circle{2}}
\put(16,10){\line(1,0){13}}
\put(30,10){\circle{2}}
\put(31,10){\line(1,0){13}}
\put(47,10){\circle*{1}}
\put(50,10){\circle*{1}}
\put(53,10){\circle*{1}}
\put(56,10){\line(1,0){13}}
\put(70,10){\circle{2}}
\put(71,9.5){\line(2,-1){13}}
\put(85,3){\circle{2}}
\put(71,10.5){\line(2,1){13}}
\put(85,17){\circle{2}}
\put(14,13){$\ell$}
\end{picture}
}
\vskip -.1cm \noindent
and $\psi_{\ell;2n}|_{U(n)} = \sum_{r+s=\ell}\psi_{r,s;n}$ where 
$\psi_{r,s;n}$ has diagram

\centerline{
\setlength{\unitlength}{.65 mm}
\begin{picture}(85,10)
\put(5,1){\circle{2}}
\put(5,4){$r$}
\put(6,1){\line(1,0){13}}
\put(20,1){\circle{2}}
\put(21,1){\line(1,0){13}}
\put(37,1){\circle*{1}}
\put(40,1){\circle*{1}}
\put(43,1){\circle*{1}}
\put(46,1){\line(1,0){13}}
\put(60,1){\circle{2}}
\put(60,4){$s$}
\put(75,0){$\times$}
\put(72,4){$s-r$}
\end{picture}.}  
\vskip .1cm \noindent
Both $\psi_{r,s;n}$ and $\psi_{r,s;n}|_{SU(n)}$ have highest weight vector
$z_1^r \bar z_1^s$.  Let $E_{r,s;n}$ denote the representation space
for $\psi_{r,s;n}$.  Now 
$\cA(G_n/K_n) = \sum_{\ell \geqq 0}\sum_{r+s=\ell}\, E_{r,s;n}$.
\medskip

{\bf Symplectic Group Cases.} 
First suppose $G_n = Sp(n) \times U(1)$.  There
are two cases: (i) the $K_n = \{1\} \times Sp(n-1) \times U(1)_{diag,n}$
and (ii) the $K_n = U(1)\times Sp(n-1)\times U(1)$.  The assertions for
case (i) will imply them for case (ii), so we may assume that
$K_n = \{1\} \times Sp(n-1) \times U(1)_{diag,n}$.  Then
$(G,K) = \varinjlim \{(G_n, K_n)\}$ has defining representation
$\pi = \varinjlim \pi_n$ where $\pi_n$ is the representation
\setlength{\unitlength}{.60 mm}
\begin{picture}(88,7)
\put(5,1){\circle{2}}
\put(5,4){{\footnotesize $1$}}
\put(6,1){\line(1,0){13}}
\put(20,1){\circle{2}}
\put(21,1){\line(1,0){13}}
\put(37,1){\circle*{1}}
\put(40,1){\circle*{1}}
\put(43,1){\circle*{1}}
\put(46,1){\line(1,0){13}}
\put(60,1){\circle{2}}
\put(61,0.5){\line(1,0){13}}
\put(61,1.5){\line(1,0){13}}
\put(75,1){\circle{2}}
\put(66,-0.2){$<$}
\put(80,0){$\times$}
\put(81,4){$1$}
\end{picture}
of $G_n$ on $\C^{2n}$.
\medskip

Note that $\pi_n$ factors through the vector representation of
$U(2n)$ on $\C^{2n}$.
We saw how $U(2n)$ acts irreducibly on the space $E_{r,s;2n}$ by the 
representation $\psi_{r,s;2n}$, which has diagram

\centerline{
\setlength{\unitlength}{.75 mm}
\begin{picture}(85,7)
\put(5,1){\circle{2}}
\put(5,4){$r$}
\put(6,1){\line(1,0){13}}
\put(20,1){\circle{2}}
\put(21,1){\line(1,0){13}}
\put(37,1){\circle*{1}}
\put(40,1){\circle*{1}}
\put(43,1){\circle*{1}}
\put(46,1){\line(1,0){13}}
\put(60,1){\circle{2}}
\put(60,4){$s$}
\put(75,0){$\times$}
\put(72,4){$s-r$}
\end{picture}.}
\vskip .1 cm \noindent  We now need two facts.  First,
$G_n\hookrightarrow U(2n)$ sends the $U(1)$ factor of $G_n$ to the center 
of $U(n)$.  Second,
$\psi_{r,s;2n}|_{Sp(n)} = \sum_{0\leqq m\leqq \min(r,s)}{'\varphi_{r,s,m;n}}$
where $'\varphi_{r,s,m;n}$ has diagram
\setlength{\unitlength}{.75 mm}
\begin{picture}(78,7)
\put(5,1){\circle{2}}
\put(-2,4){{\footnotesize $r+s-2m$}}
\put(6,1){\line(1,0){13}}
\put(20,1){\circle{2}}
\put(21,4){{\footnotesize $m$}}
\put(21,1){\line(1,0){13}}
\put(37,1){\circle*{1}}
\put(40,1){\circle*{1}}
\put(43,1){\circle*{1}}
\put(46,1){\line(1,0){13}}
\put(60,1){\circle{2}}
\put(61,0.5){\line(1,0){13}}
\put(61,1.5){\line(1,0){13}}
\put(75,1){\circle{2}}
\put(66,-0.2){$<$}
\end{picture}. 
That gives us $\psi_{r,s;2n}|_{Sp(n)U(1)} = \sum_{0 \leqq m \leqq \min(r,s)}
\varphi_{r,s,m;n}$ where $\varphi_{r,s,m;n}$ is the representation of 
$Sp(n)U(1)$ with diagram 
\setlength{\unitlength}{.75 mm}
\begin{picture}(92,7)
\put(5,1){\circle{2}}
\put(-2,4){{\footnotesize $r+s-2m$}}
\put(6,1){\line(1,0){13}}
\put(20,1){\circle{2}}
\put(21,4){{\footnotesize $m$}}
\put(21,1){\line(1,0){13}}
\put(37,1){\circle*{1}}
\put(40,1){\circle*{1}}
\put(43,1){\circle*{1}}
\put(46,1){\line(1,0){13}}
\put(60,1){\circle{2}}
\put(61,0.5){\line(1,0){13}}
\put(61,1.5){\line(1,0){13}}
\put(75,1){\circle{2}}
\put(66,-0.2){$<$}
\put(85,0){$\times$}
\put(82,4){$s-r$}
\end{picture},
and $'\varphi_{r,s,m;n}$ has the same representation space 
(call it $E_{r,s,m;n}$) as $\varphi_{r,s,m;n}$.  The
$E_{r,s,m;n}$ are irreducible and inequivalent under $Sp(n)U(1)$;
in other words the irreducible representations $\varphi_{r,s,m;n}$
all are mutually inequivalent.
Note, however, that $'\varphi_{r,s,m;n} \simeq {'\varphi_{r+t,s-t,m;n}}$
for all $t$ such that $r+t, s-t \geqq 0$; this reflects the fact that
$(Sp(n),Sp(n-1))$ is not a Gelfand pair.
\medskip

To trace the inclusions let $\{z_1, \dots , z_{2n}\}$ be the coordinates
of $\C^{2n}$, all weight vectors, where $z_1$ is the highest weight vector, 
$z_2 = e_{-\alpha_1}z_1$ is the next highest, and so on, and the antisymmetric
bilinear invariant of $Sp(n)$ on $\C^{2n}$ is
$v_n(z,w) = \sum_1^n (z_{2i-1} w_{2i} - z_{2i} w_{2i-1})$. 
Then $z_1^\ell$ is the highest weight vector of
\setlength{\unitlength}{.60 mm}
\begin{picture}(80,7)
\put(5,1){\circle{2}}
\put(5,4){{\footnotesize $\ell$}}
\put(6,1){\line(1,0){13}}
\put(20,1){\circle{2}}
\put(21,1){\line(1,0){13}}
\put(37,1){\circle*{1}}
\put(40,1){\circle*{1}}
\put(43,1){\circle*{1}}
\put(46,1){\line(1,0){13}}
\put(60,1){\circle{2}}
\put(61,0.5){\line(1,0){13}}
\put(61,1.5){\line(1,0){13}}
\put(75,1){\circle{2}}
\put(66,-0.2){$<$}
\end{picture} and
$\Lambda^2 \C^{2n}$ is the sum $\Lambda^2_0 \C^{2n} \oplus v_n\C$ of its
irreducible component and its trivial component 
under the action of $Sp(n)$.  Here $v_n$ has matrix
$\diag \left \{ 
\left [\begin{smallmatrix} 0 & 1 \\ -1 & 0 \end{smallmatrix} \right ]\dots
\left [\begin{smallmatrix} 0 & 1 \\ -1 & 0 \end{smallmatrix} \right ]\right \}$ 
and we work with the maximal toral subalgebra that consists of all matrices
$\diag \left \{ a_1, -a_1; \dots ; a_n,-a_n \right \}$; thus
the highest weight vector on $\Lambda^2_0 \C^{2n}$
is $s_n(z,w) = z_1w_3 - z_3w_1$.  Now $s_n^m$ is the highest weight vector of
\setlength{\unitlength}{.60 mm}
\begin{picture}(80,7)
\put(5,1){\circle{2}}
\put(6,1){\line(1,0){13}}
\put(20,1){\circle{2}}
\put(20,4){{\footnotesize $m$}}
\put(21,1){\line(1,0){13}}
\put(37,1){\circle*{1}}
\put(40,1){\circle*{1}}
\put(43,1){\circle*{1}}
\put(46,1){\line(1,0){13}}
\put(60,1){\circle{2}}
\put(61,0.5){\line(1,0){13}}
\put(61,1.5){\line(1,0){13}}
\put(75,1){\circle{2}}
\put(66,-0.2){$<$}
\end{picture}.
The corresponding highest weight vector of $\varphi_{r,s,m;n}$ is
$z_1^{r-m} \bar z_1^{s-m} s_n^m$.  Now the restriction 
$\cA(G_{n+1}/K_{n+1}) \to \cA(G_n/K_n)$ maps the highest weight vector
$z_1^{r-m} \bar z_1^{s-m} s_{n+1}^m$ of $\varphi_{r,s,m;n+1}$ to the highest 
weight vector $z_1^{r-m} \bar z_1^{s-m} s_n^m$.  This proves that the
representation of $G$
on $\cA(G/K)$ is limit--aligned.  Theorem \ref{reduction-cor} shows
that $\varinjlim \cA(G_n/K_n)$ is the multiplicity free direct sum of
lim--irreducible $G$--modules $E_{r,s,m} := \varinjlim E_{r,s,m;n}$.
\medskip

Finally, we suppose $G_n = Sp(n)\times Sp(1)$.  Again there are two cases:
(i) the $K_n = \{1\} \times Sp(n-1) \times Sp(1)_{diag,n}$
and (ii) the $K_n = U(1)\times Sp(n-1)\times Sp(1)$.  The function algebras
and group actions in case (ii) are exactly the same as those of the setting
$(G_n,K_n) = (Sp(n)\times U(1), U(1)\times Sp(n-1)\times U(1))$ above, where
the assertions are proved.  Thus we need only consider case (i),
$K_n = \{1\} \times Sp(n-1) \times Sp(1)_{diag,n}$.  Then
$(G,K) = \varinjlim \{(G_n, K_n)\}$ has defining representation
$\pi = \varinjlim \pi_n$ described in ``Symplectic 1'' above.  Those $\pi_n$
satisfy the condition of Theorem \ref{reduction-cor} because
$Sp(n)\times Sp(1)$ simply puts together representation spaces
$E_{r-m,s-m,m;n}$ of $Sp(n)\times U(1)$ on $\cA(Sp(n)U(1)/Sp(n-1)U(1))$.  This 
assembly maintains total degree $\ell = (r-m)+(s-m)+2m$, views the
$U(1)$ factor of $Sp(n)\times U(1)$ as a maximal torus of the $Sp(1)$ factor of
$Sp(n)\times Sp(1)$, and sums the spaces for the 
\setlength{\unitlength}{.60 mm}
\begin{picture}(10,7)
\put(2,1){$\times$}
\put(0,4){\footnotesize $s-r$}
\end{picture}
to form the space for the irreducible representation (call it $\beta_\ell$) 
of $Sp(1)$ of degree $\ell +1$.  It has diagram
\setlength{\unitlength}{.60 mm}
\begin{picture}(4,7)
\put(2,1){\circle{2}}
\put(0,3){\footnotesize $\ell$}
\end{picture}.
\, Now the irreducible spaces for $Sp(n)\times Sp(1)$ are the
$F_{\ell,m;n}:= \sum_{r+s=\ell} E_{r-m,s-m,m;n}$ and the corresponding
representations are the
$\varphi_{\ell,m,n} := \sum_{r+s=\ell} {'\varphi_{r-m,s-m,m;n}}$.
This proves that the representation of $G$
on $\cA(G/K)$ is limit--aligned.  Theorem \ref{reduction-cor} shows
that $\varinjlim \cA(G_n/K_n)$ is the multiplicity free direct sum of
lim--irreducible $G$--modules $F_{\ell,m} := \varinjlim F_{\ell,m;n}$.
\medskip

We have proved Theorem \ref{spec}.
\smartqed\qed\end{proof}

\begin{rem}\label{additional} {\rm
Alternatively the systems (d), (e) and (f) from the list (\ref{kramer-seq}),
and also (a) when the $\{p_n\}$ are bounded, can be treated by the method 
of Sections \ref{sec6} and \ref{sec7}.  That gives an alternative
proof of the multiplicity free property for the pairs
$$
\begin{aligned}
&{\rm (1)} \phantom{X}
SU(p+\infty)/(SU(p)\times SU(\infty)) \text{ for } 1 \leqq p \leqq \infty, \\
&{\rm (2)} \phantom{X} SU(1+2\infty)/(U(1)\times Sp(\infty)), \\
&{\rm (3)} \phantom{X} SO(1 + 2\infty)/U(\infty), \text{ and } \\
&{\rm (4)} \phantom{X} Sp(1+\infty)/(U(1)\times Sp(\infty))
\end{aligned}
$$
of Theorem \ref{cor-gf1-mfree}.
} \hfill $\diamondsuit$
\end{rem}

\section{Summary}\label{sec9}
We have proved that the regular representations of $G$ on $\cA(G/K)$, 
$\cC(G/K)$ and $L^2(G/K)$ are multiplicity free discrete direct sums of
lim--irreducible representations in the following cases.  In addition, 
in those cases it is always permissible to enlarge the groups $K_n$,
say to $F\cdot K_n$ where $F$ is a closed subgroup of the normalizer
$N_{G_n}(K_n)$, because $\cA(G_n/[F\cdot K_n])$ is a $G_n$--submodule
of $\cA(G_n/K_n)$.
\medskip

\underline{Limits of riemannian symmetric spaces.}  We have the 
multiplicity free property for the thirteen cases 
described in Theorem \ref{cor-symm-mfree}, as well as some obvious variations.
The latter include
$$
\begin{aligned}
SO(\infty)\times SO(\infty)/\text{diag\,} SO(\infty) 
	= &\varinjlim SO(n)\times SO(n)/\text{diag\,} SO(n)\ \text{ and } \\
SO(p+\infty)/[S(O(p) \times O(\infty))]
	= &\varinjlim SO(n)/[S(O(p) \times O(n-p))]\ .
\end{aligned}
$$

\underline{Limits of a few systems of Gelfand pairs.}  We have the 
multiplicity free property for the four cases described in 
Theorem \ref{cor-gf1-mfree}, 
$$
\begin{aligned}
&{\rm (1)} \phantom{X}
SU(p+\infty)/(SU(p)\times SU(\infty)) \text{ for } 1 \leqq p \leqq \infty, \\
&{\rm (2)} \phantom{X} SU(1+2\infty)/(U(1)\times Sp(\infty)), \\
&{\rm (3)} \phantom{X} SO(1 + 2\infty)/U(\infty), \text{ and } \\
&{\rm (4)} \phantom{X} Sp(1+\infty)/(U(1)\times Sp(\infty)).
\end{aligned}
$$
We also have the multiplicity free property for spaces that 
interpolate between 
$$
(SU(p+\infty), SU(p)\times SU(\infty))
$$ 
and the limit Grassmannian
$$
(U(p+\infty), \varinjlim U(p)\times U(n).
$$  
Fix a closed subgroup $F$ 
of $U(1)$.  Then we have the multiplicity free property for the 
pairs $(G,K) = \varinjlim \{(G_n,K_n)\}$ where $G_n = SU(p+n)$
and 
$$
K_n = \left \{ \left . \left ( \begin{smallmatrix} k_n' & 0 \\ 0 & k_n'' 
\end{smallmatrix} \right ) \right |
k_n' \in U(p), k_n'' \in SU(n), \det k_n' \in F \right \}.
$$

\underline{Limits of Gelfand pairs related to spheres and
Grassmann manifolds.}  We have the multiplicity free property for the 
six cases described in Theorem \ref{spec},
four of which are nonsymmetric, as well as some obvious variations.  Fix a
closed subgroup $F$ of $U(1)$; it can be any finite cyclic group or
the entire circle group $U(1)$.  As a result we have the multiplicity free
property for the nonsymmetric pairs
$$
\begin{aligned}
&SU(\infty)/[SU(p)\times SU(\infty - p)]\\
	&\phantom{XXXX}= \varinjlim SU(n)/[SU(p)\times SU(n - p)]\ ,\\
&[Sp(\infty)\times U(1)]/[F \times Sp(\infty - 1) \times U(1)_{diag}]\\
	&\phantom{XXXX}= 
	  \varinjlim [Sp(n)\times U(1)]/[F\times Sp(n-1)\times U(1)_{diag}]
		\ ,\\
&[Sp(\infty)\times Sp(1)]/[\{1\} \times Sp(\infty - 1) \times Sp(1)_{diag}]\\
	&\phantom{XXXX}= \varinjlim [Sp(n)\times Sp(1)]/
	[\{1\}\times Sp(n-1)\times Sp(1)_{diag}] \ , \text{ and}\\
&[Sp(\infty)\times Sp(1)]/[\{\pm1\}\times Sp(\infty -1)\times Sp(1)_{diag}]\\
	&\phantom{XXXX}= \varinjlim [Sp(n)\times Sp(1)]/
	[\{\pm1\}\times Sp(n-1)\times Sp(1)_{diag}] \ .
\end{aligned}
$$
\medskip

\underline{What we don't have.}
There is a huge number of direct systems $\{(G_n,K_n)\}$ 
of Gelfand pairs where the $G_n$ are compact connected Lie groups.
We have only verified the multiplicity free condition for a few of
them.  We have not, for example, checked it for the interesting cases
$$
G_n = SU(2n+1) \text{ and } K_n = F \times Sp(n), \ \ F \subset U(1) 
\text{ finite cyclic,}
$$
and
$$
G_n = SO(2n) \text{ and } K_n = F \times SU(n), \ n \text{ odd, } n \geqq 3.
$$
Also, we have not checked it for the {\sl very} interesting case
$$
G_n = Sp(a_n) \times Sp(b_n) \text{ and }
K_n = Sp(a_n -1)\times Sp(1)\times Sp(b_n -1),
$$
which is a prototype for nonsymmetric irreducible direct systems
$\{(G_n,K_n)\}$ with the $G_n$ semisimple but not simple.
In that case $K_n \hookrightarrow G_n$ is given by
$(k_1,a,k_2) \mapsto \left ( 
\left ( \begin{smallmatrix} k_1 & 0 \\ 0 & a \end{smallmatrix} \right ),
\left ( \begin{smallmatrix} a & 0 \\ 0 & k_2 \end{smallmatrix} \right )
\right )$, so $G_n/K_n$ fibers over $P^{a_n -1}(\H)\times P^{b_n -1}(\H)$
with fiber $(Sp(1)\times Sp(1))/(diagonal) = S^3$.
\vfill\pagebreak

\printindex


\begin{thebibliography}{XXXXXXI}

\bibitem[DPW2002]{DPW2002}
I. Dimitrov, I. Penkov \& J. A. Wolf,
A Bott--Borel--Weil theory for direct limits of algebraic groups,
Amer. J of Math. {\bf 124} (2002), 955--998.

\bibitem[Fa2006]{Fa2006}
J. Faraut,
Infinite dimensional harmonic analysis and probability,
in ``Probability Measures on Groups: Recent Directions and Trends,''
ed. S. G. Dani \& P. Graczyk, Narosa, New Delhi, 2006.

\bibitem[Kr1979]{Kr1979}
M. Kr\" amer,
Sph\" arische Untergruppen in kompakten zusammenh\" angenden Liegruppen,
Compositio Math. {\bf 38} (1979), 129--153.

\bibitem[Ol1990]{Ol1990}
G. I. Ol'shanskii,
Unitary representations of infinite dimensional pairs $(G,K)$ and the
formalism of R. Howe, in ``Representations of Lie Groups and Related
Topics, ed. A. M. Vershik \& D. P. Zhelobenko,'' Advanced Studies
Contemp. Math. {\bf 7}, Gordon \& Breach, 1990.

\bibitem[W2005]{W2005}
J. A. Wolf,
Direct limits of principal series representations,  Compositio
Mathematica, {\bf 141} (2005), 1504--1530.

\bibitem[W2007]{W2007}
J. A. Wolf,
Harmonic Analysis on Commutative Spaces, 
Math. Surveys \& Monographs vol. 142, Amer. Math. Soc., 2007.

\bibitem[W2008]{W2008}
J. A. Wolf,
Infinite dimensional multiplicity free spaces II:
Limits of commutative nilmanifolds, to appear.

\bibitem[Ya2004]{Ya2004}
O. S. Yakimova,
Weakly symmetric riemannian manifolds with reductive isometry group,
Math. USSR Sbornik {\bf 195} (2004), 599--614.

\end{thebibliography}
\end{document}